\documentclass[10pt]{amsart}
\usepackage{latexsym}
\usepackage{amsfonts}

\usepackage{amsmath,amsthm,amssymb}

\newtheorem{theor}{Theorem}

\newtheorem{corol}[theor]{Corollary}

\newtheorem{thm}{Theorem}[section]
\newtheorem{lem}[thm]{Lemma}
\newtheorem{cor}[thm]{Corollary}

\newtheorem{prop}[thm]{Proposition}
\theoremstyle{definition}
\newtheorem{defn}[thm]{Definition}

\newtheorem{conv}[thm]{Convention}
\newtheorem{rem}[thm]{Remark}
\newtheorem{exmp}[thm]{Example}

\newtheorem{prob}[thm]{Problem}

\newcommand\m{\mu_A}

\begin{document}

\title[Generic stretching factors]{The Subadditive Ergodic Theorem and
  generic stretching factors for free group automorphisms}

\author[V.~Kaimanovich]{Vadim Kaimanovich}

\address{\tt IRMAR, Universit\'{e} Rennes-1, Campus de Beaulieu, 35042, Rennes Cedex,
France;} \email{\tt kaimanov@univ-rennes1.fr}

\author[I.~Kapovich]{Ilya Kapovich}

\address{\tt Department of Mathematics, University of Illinois at
  Urbana-Champaign, 1409 West Green Street, Urbana, IL 61801, USA
  \newline http://www.math.uiuc.edu/\~{}kapovich/} \email{\tt
  kapovich@math.uiuc.edu}

\author[P.~Schupp]{Paul Schupp}

\address{\tt Department of Mathematics, University of Illinois at
  Urbana-Champaign, 1409 West Green Street, Urbana, IL 61801, USA
  \newline http://www.math.uiuc.edu/People/schupp.html}
\email{schupp@math.uiuc.edu}

\thanks{The second and the third author were supported by the NSF
  grant DMS\#0404991 and the NSA grant DMA\#H98230-04-1-0115}

\begin{abstract}
  Given a free group $F_k$ of rank $k\ge 2$ with a fixed set of free
  generators we associate to any homomorphism $\phi$ from $F_k$ to a
  group $G$ with a left-invariant semi-norm a generic stretching
  factor, $\lambda(\phi)$, which is a non-commutative generalization
  of the translation number.  We concentrate on the situation when
  $\phi:F_k\to Aut(X)$ corresponds to a free action of $F_k$ on a
  simplicial tree $X$, in particular, when $\phi$ corresponds to the
  action of $F_k$ on its Cayley graph via an automorphism of $F_k$. In
  this case we are able to obtain some detailed ``arithmetic''
  information about the possible values of $\lambda=\lambda(\phi)$. We
  show that $\lambda \ge 1$ and is a rational number with
  $2k\lambda\in \mathbb Z\left[ \frac{1}{2k-1} \right]$ for every
  $\phi\in Aut(F_k)$.  We also prove that the set of all
  $\lambda(\phi)$, where $\phi$ varies over $Aut(F_k)$, has a gap
  between $1$ and $1+\frac{2k-3}{2k^2-k}$, and the value $1$ is
  attained only for ``trivial'' reasons. Furthermore, there is an
  algorithm which, when given $\phi$, calculates $\lambda( \phi)$.
\end{abstract}


\maketitle

\tableofcontents

\section{Introduction}\label{intro}

\subsection{Random subgroup distortion and growth of random automorphisms}

Let $G$ be a finitely generated group with a word metric $d_S$
determined by a finite generating set $S$ and write $|g|_S:=d_S(1,g)$
for $g\in G$.  Recall that if $H=\langle T\rangle$ is a subgroup of
$G$ generated by a finite set $T$, then a function $f$ is said to be a
\emph{distortion function} of $H$ in $G$ if for every $h\in H$ we have
$|h|_T\le f(|h|_S)$. The subgroup $H$ is \emph{quasi-isometrically
  embedded} in $G$ if and only if for some (and hence for all) choices
of $S,T$ there is a linear distortion function for $H$ in $G$, that is
if the ratio $\frac{|h|_T}{|h|_S}$ is bounded on $H\setminus\{1\}$.

The \emph{translation number} of an element $g\in G$ is defined as
\[
\lambda(g) = \lambda_S(g) =\lim_{n\to\infty} \frac{|g^n|_S}{n}
\]
and the limit exists by the subadditivity of the sequence $|g^n|_S$.
If $g$ has infinite order, then the cyclic subgroup $\langle g\rangle$
is quasi-isometrically embedded in $G$ if and only if $\lambda_S(g)>0$
for some (and hence for any) finite generating set $S$ of $G$.

The study of ``random'' or ``generic'' behavior is currently an
increasingly active area of investigation in many different
group-theoretic contexts. (See, for example,
\cite{Grom1,Ol92,Grom2,Ch94,Ch95,Ch00,Che96,Che98,A1,A2,A3,AO,Z,KMSS,KMSS1,KS,Gh,Oliv}
).  In this paper we concentrate on algebraic and geometric
consequences of subadditivity, specifically of Kingman's Subadditive
Ergodic Theorem.

We investigate a ``noncommutative analogue'' of translation number
which is defined for a ``mapped-in'' subgroup $H=\phi(F)$, where
$\phi:F\to G$ is a homomorphism of a free nonabelian group $F=F(A)$ of
finite rank into a group $G$ with generating set $S$.  Namely, there
is a number $\lambda=\lambda(\phi,A,S)\ge 0$ such that for long
``random'' freely reduced words $w\in F(A)$ we have
$\frac{|\phi(w)|_S}{|w|_A}\approx \lambda$.  (Instead of the word
metric $d_S$ one could actually take an arbitrary semi-norm on $G$.)

\medskip

Throughout this paper we fix the notation that $F=F(A)$ is the free
group with basis $A=\{a_1,\dots, a_k\}$ where $k\ge 2$. For any $w\in
F$ let $|w|$ denote the length of the unique freely reduced word over
$A^{\pm 1}$ representing $w$. We identify the hyperbolic boundary
$\partial F$ with the set of all \emph{geodesic rays} from $1\in F$ in
the Cayley graph $\Gamma(F,A)$ of $F$, that is, $\partial F$ is the
set of all semi-infinite freely reduced words over $A^{\pm 1}$ endowed
with the standard topology. The space $\partial F$ can be identified
with the \emph{space of ends} or the \emph{hyperbolic boundary} of
$F$.

The Borel $\sigma$-algebra $\mathcal F$ on $\partial F$ is generated
by the \emph{cylinder sets} $Cyl_A(v),\,v\in F$, where $Cyl_A(v)$
consists of all infinite rays $\omega\in \partial F$ that begin with
$v$.  The \emph{uniform} Borel probability measure $\m$ on $\partial
F$ corresponding to $A$ is defined by assigning equal weights to all
cylinders based on the words on the same length. That is,
$$
\m(Cyl_A(v))=\frac{1}{2k(2k-1)^{|v|-1}} \qquad\forall\, v\in
F\setminus\{1\} .
$$
Note that although the boundary $\partial F$ could be defined
without referring to a particular generating set $A$, the uniformity
of the measure $\m$ \emph{does} depend on the choice of $A$. The fact
that the uniform measures corresponding to two different free
generating sets may well be singular respect to each other is actually
the cornerstone of our approach.  (See \cite{Fur} for a detailed
discussion of this phenomenon in the general context of
word-hyperbolic groups and of the Patterson-Sullivan measures
corresponding to geometric actions of such groups on Gromov-hyperbolic
spaces.)

A ray $\omega\in\partial F$ can be thought of as a
\emph{non-backtracking edge-path} in $\Gamma(F,A)$ starting from the
identity $1$. We denote by $\omega_n$ the vertex on $\omega$ at
distance $n$ from $1$. The measure space $(\partial F,\m)$ then
becomes the \emph{space of sample paths} of the \emph{non-backtracking
  simple random walk} (NBSRW) on the Cayley graph of $F$.  This is the
Markov chain on $F$ whose transition probabilities $\pi_w,\,w\in F$
are equidistributed among the neighbors of $w$ which are strictly
further from the group identity. By choosing a random $\m$-distributed
point $\omega\in\partial F$ we may think about $\omega_n$ as a
``random'' (with respect to the NBSRW) freely reduced word of length
$n$ in $F$.

%\newpage

In Section~2 we prove:

\begin{theor}\label{thm:stretch}
  Let $F=F(a_1,\dots, a_k)$ with $k\ge 1$, and let $\m$ be the uniform
  Borel probability measure on $\partial F$ corresponding to the basis
  $A=\{a_1,\dots, a_k\}$.

  Let $\phi:F\to G$ be a homomorphism to a group $G$ endowed with a
  semi-norm, that is, a nonnegative function $|\cdot|_G$ on $G$
  satisfying $|gh|_G\le |g|_G+|h|_G$ for all $g,h\in H$.

  Then:

\begin{enumerate}
\item There exists a real number $\lambda\ge 0$ such that
\[
\lim_{n\to\infty} \frac{|\phi(\omega_n)|_G}{n}=\lambda
\]
$\m$-a.e. and in $L^1(\partial F,\m)$.
\item Suppose further that the image group $\phi(F)$ is non-amenable,
  and that the sequence $b_n=\#\{ g\in \phi(F): |g|_G\le n\}$ grows at
  most exponentially. Then $\lambda>0$.
\end{enumerate}

\end{theor}
Note that the requirement of at most exponential growth of the $b_n$
is automatically satisfied if the group $G$ is finitely generated, and
$|\cdot|_G$ is the word metric on $G$ determined by a finite
generating set.

Theorem~\ref{thm:stretch} says that for a long ``random'' freely
reduced element $w\in F$ we have $\frac{|\phi(w)|_G}{|w|}\approx
\lambda$. For this reason we shall call the number
$\lambda=\lambda(\phi, A,|\cdot|_G)$, whose existence is provided by
part (1) of Theorem~\ref{thm:stretch}, the \emph{generic stretching
  factor} of $\phi$ with respect to the free basis $A$ of $F$ and the
semi-norm $|\cdot|_G$.

We deduce Theorem~\ref{thm:stretch} from the fact that the sample
paths of the usual simple random walk on the group $F$ asymptotically
follow geodesics and the well-known results on the linear rate of
escape of random walks on groups \cite{KV}, \cite{K}. We also give an
alternative direct argument proof of part (1) of
Theorem~\ref{thm:stretch} applying Kingman's Subadditive Ergodic
Theorem \cite{King}.  Part (2) of Theorem~\ref{thm:stretch} can also
be proved using the results of Cohen~\cite{Cohen},
Grigorchuck~\cite{Gri} and Woess~\cite{Woess83} on co-growth in
groups.

\medskip

\begin{exmp}[Stretching factors for isometric actions]\label{exmp:metric}
  A typical example of a semi-norm $|\cdot|_G$ comes from isometric
  group actions on metric spaces. Namely, let $X$ be a metric space
  with basepoint $x\in X$.  For an isometry $g$ of $X$ define
  $|g|_x:=d(x,gx)$. The triangle inequality implies that
  $|g_1g_2|_x\le |g_1|_x+|g_2|_x$, so that $|\cdot|_x$ is a semi-norm
  on $G=Isom(X)$. Suppose $F=F(a_1,\dots, a_k)$ acts by isometries on
  $X$ by a homomorphism $\phi: F\to G$. It is easy to see that in this
  case $\lambda(\phi, A, |\cdot|_x)$ does not depend on the choice of
  a base-point $x\in X$ and is determined by the action $\phi$ and the
  choice of the basis $A$ of $F$. In this case we shall denote
  $\lambda(\phi, A, |\cdot|_x)$ by $\lambda(\phi,A)$, or just by
  $\lambda(\phi)$ if the choice of $A$ is fixed.
\end{exmp}

\begin{exmp}[Random Subgroup Distortion]
  Let $H\le G$ be finitely generated groups with finite generating
  sets $A\subset H$ and $S\subset G$, respectively. Denote the
  associated length functions by $|\cdot|_G$ and $|\cdot|_H$. Now $H$
  is a quotient of $F=F(A)$.  Let $\phi:F(A)\to G$ be composition of
  this quotient map with the inclusion of $H$ into $G$. Then
  $|\phi(w)|_H\le |w|$ for any $w\in F$. If the group $H$ is
  non-amenable then by Theorem~\ref{thm:stretch} for a long ``random''
  freely reduced word $w$ in $F(A)$

\[
\frac{|\phi(w)|_H}{|w|}\approx \lambda_1>0, \quad \text{ and } \quad
\frac{|\phi(w)|_G}{|w|}\approx \lambda_2>0,
\]
and therefore
\[
\frac{|\phi(w)|_H}{|\phi(w)|_G}\approx \frac{\lambda_1}{\lambda_2},
\]
where the constants $\lambda_1,\lambda_2>0$ do not depend on $w$.
Thus, informally speaking, Theorem~\ref{thm:stretch} implies that any
nonamenable finitely generated subgroup $H$ of a finitely generated
group $G$ has ``generic-case'' linear distortion in $G$.
\end{exmp}

\begin{exmp}[Normal Forms]
  Let $G$ be a nonamenable group with a finite generating set $A$ and
  the associated length function $|\cdot|_G$. We will denote by
  $\overline w$ the element of $G$ represented by a word $w$ in the
  alphabet $A\cup A^{-1}$.

  Let $L\subseteq (A\cup A^{-1})^{\ast}$ be a set of \emph{normal
    forms} (not necessarily unique) for elements of $G$, that is
  $\overline L=G$. Consider, for instance, the case where $L$ is an
  automatic language for $G$. By Theorem~\ref{thm:stretch} there is
  $\lambda>0$ such that for a random long freely reduced word $w\in
  F(A)$ we have $\frac{|\overline{w}|_G}{|w|}\approx \lambda$. Let
  $w_L\in L$ be a word representing the same element of $G$ as $w$.
  Then
\[
|w_L|\ge |\overline{w}|_G
\]
and hence
\[
\frac{|w_L|}{|w|}\ge \frac{|\overline{w}|_G}{|w|}\approx \lambda>0.
\]

Thus for a long random word $w\in F(S)$ when we take $\overline w$ to
a normal form $w_L\in L$, the ratio $\frac{|w_L|}{|w|}$ is separated
from zero. This conclusion applies to a number of experimental
observations, such as those obtained by Dehornoy~\cite{Deh} in the
case of braid groups.
\end{exmp}

Theorem~\ref{thm:stretch} has implications regarding the growth of
\emph{random automorphisms}. Let $G$ be a finitely generated group
with a fixed word metric corresponding to a finite generating set $S$.
Let $\phi\in Aut(G)$ be an automorphism. We define the \emph{norm} of
$\phi$ with respect to $S$ as
\[
||\phi||=||\phi||_S:=\max_{s\in S} |\phi(s)|_S.
\]

Then for any $g\in G$ we have $|\phi(g)|_S\le ||\phi||_S |g|_S$ and
hence
\[
||\phi||=\sup_{g\in G, g\ne 1} \frac{|\phi(g)|_S}{|g|_S}.
\]

For an individual $\phi$ the sequence $\log ||\phi^n||$ is subadditive
and therefore the following limit (sometimes called the \emph{growth
  entropy of $\phi$}) exists:
\[
\nu(\phi):=\lim_{n\to\infty} \frac{\log ||\phi^n||}{n}.
\]

It turns out that this notion has a generalization for an arbitrary
finitely generated subgroup of $Aut(G)$:

\begin{theor}\label{aut}
  Let $G$ be a nontrivial finitely generated group with a word-metric
  $d_S$ corresponding to a finite generating set $S$.  Let $H\le
  Aut(G)$ be a noncyclic subgroup with a finite generating set $T$.
  Then:

\begin{enumerate}
\item There is $\nu=\nu(H)=\nu(H,T,S)\ge 0$ such that for a
  non-backtracking simple random walk $\phi_n$ on the Cayley graph of
  $H$ with respect to $T$ we have

\[
\lim_{n\to\infty}\frac{\log ||\phi_n||_S}{n}=\nu
\]
almost surely and in $L^1$.

\item If $G$ has polynomial growth and $H$ is non-amenable then
  $\nu(H,T,S)>0$.
\end{enumerate}
\end{theor}

Note that $\nu(H,T,S)>0$ means that $||\phi_n||_S$ grows exponentially
with $n$ for a "random" automorphism $\phi_n$.

\begin{corol}\label{cor:aut}
  Let $F$ be a free group of finite rank $k\ge 2$ and let $H\le
  Aut(F)$ be a finitely generated group of automorphisms of $F$ such
  that the image $H'$ of $H$ in $Aut(F_{ab})\cong GL(k,{\mathbb Z})$
  is non-amenable. Then for any finite generating set $S$ of $F$ and
  for any finite generating set $T$ of $H$ we have $\nu(H,T,S)>0$.
\end{corol}

By the Tits alternative a subgroup of $GL(k,{\mathbb Z})$ is either
virtually solvable (and hence amenable) or it contains a free subgroup
of rank two (and hence is nonamenable). Thus in the above corollary we
could replace the assumption that $H'$ is nonamenable by the
requirement that $H'$ is not virtually solvable.

\subsection{Free actions on trees:  Two interpretation of stretching factors}

In the context of free and discrete isometric actions of free groups
on $\mathbb R$-trees (cf. Example~\ref{exmp:metric}).  generic
stretching factors are related to Bonahon's notion~\cite{Bo86,Bo88} of
the intersection number between geodesic currents on hyperbolic
surfaces.  If $G$ is a non-elementary word-hyperbolic group, a
\emph{geodesic current} on $G$ is a $G$-invariant positive Borel
measure on $\partial^2G:=\{(x,y)|x,y\in \partial G, x\ne y\}$. The
space of all geodesic currents on $G$, endowed with the
weak-$\ast$-topology, is denoted by $Curr(G)$.  (See \cite{Bo91,Ka1}
for a detailed discussion on the subject.)

Every nontrivial conjugacy class $[g]$ in $G$ defines an associated
``counting'' current $\eta_{[g]}$ on $G$. When $S$ is a closed surface
of negative Euler characteristic and $G=\pi_1(S)$, Bonahon proved that
the notion of geometric intersection number between free homotopy
classes of essential closed curves on $S$ (that is, between nontrivial
conjugacy classes of $G$) extends to a bilinear continuous
``intersection form'' \[i:Curr(G)\times Curr(G)\to \mathbb R.\] Note
that in this case $\partial G=\partial \mathbb H^2=\mathbb S^1$. For
every hyperbolic structure $\rho$ on $S$ there is an associated
\emph{Liouville current} $L_{\rho}\in Curr(G)$ (see \cite{Bo86}).
Bonahon's construction has the following natural property: if $\rho$
is as above and $[g]$ is a nontrivial conjugacy class in $G$ then
$i(L_{\rho}, \eta_{[g]})=\ell_{\rho}(g)$. Here $\ell_{\rho}: G\to
\mathbb R$ is the \emph{length spectrum} of $\rho$. Thus
$\ell_{\rho}(g)$ is equal to the translation length of $g$ as an
isometry of $\mathbb H^2=\widetilde{(S,\rho)}$ and it is also equal to
the $\rho$-length of the shortest curve of the free homotopy class of
closed curves on $S$ corresponding to $[g]$. It turns out that the
intersection number $i(L_\rho,L_{\rho'})$ between Liouville currents
corresponding to two hyperbolic structures $\rho,\rho'$ can be
interpreted as the generic stretching factor of a long random closed
geodesic on $(S,\rho)$ with respect to $\rho'$. Namely, let $p\in S$
and let $v$ be a random unit tangent vector at $p$ on $(S,\rho)$. For
every $n\ge 1$ let $\alpha_n$ be the geodesic of length $n$ on
$(S,\rho)$ with origin $p$ and with the tangent vector $v$ at $p$.
Let $\beta_n$ be a geodesic from the terminus of $\alpha_n$ to $p$ of
length $\le Diam(S,\rho)$. Then $\gamma_n=\alpha_n\beta_n$ is a closed
curve on $S$. Bonahon's results imply that
\[
\lim_{n\to\infty}
\frac{\ell_{\rho'}([\gamma_n])}{\ell_{\rho}([\gamma_n])}=\lim_{n\to\infty}
\frac{\ell_{\rho'}([\gamma_n)])}{n}=\frac{i(\rho,\rho')}{\pi^2|\chi(S)|}.
\]

It turns out that a version of this interpretation applies in the
context of free groups acting on trees. Let $F$ be a free group of
finite rank $k\ge 2$ In~\cite{Ka,Ka1} Kapovich investigated a natural
``intersection form'' $I: FLen(F)\times Curr(F)\to \mathbb F$, where
$FLen(F)$ is the space of hyperbolic length functions corresponding to
free and discrete isometric actions of $F$ on $\mathbb R$-trees. This
form still has the natural property that for any nontrivial conjugacy
class $[g]$ in $F$ and any $\ell\in FLen(F)$ we have $I(\ell,
\eta_{[g]}=\ell(g)$.  Let $A$ be a free basis of $F$ and let $\ell\in
FLen(F)$ be realized by a free and discrete isometric action
$\phi:F\to Isom(X)$ of $F$ on an $\mathbb R$-tree $X$. Let $\mu_A$ be
the uniform measure on $\partial F$ corresponding to $A$. The measure
$\mu_A$ on $\partial F$ determines a \emph{uniform} current $\nu_A\in
Curr(F)$ that is analogous to the Liouville current corresponding to a
hyperbolic structure on a surface.  As shown in \cite{Ka1}, similarly
to Bonahon's situation, we have
\[
I(\ell, \nu_A)=\lambda_A(\phi).
\]

Generic stretching factors are also related to the notion of the
Hausdorff dimension of a measure with respect to a metric.  If $\mu$
is a measure on a metric space $(M,d)$, the \emph{Hausdorff dimension
  of $\mu$ with respect to $d$}, denoted ${\bf HD}_d(\mu)$ (or just
${\bf HD}(\mu)$), is defined as the infimum of Hausdorff dimensions of
subsets of $(M,d)$ of full measure $\mu$.

In \cite{Kaim98} Kaimanovich proved that for the harmonic measure
$\nu$ on $\partial T$ associated to a regular Markov operator $P$ with
a positive rate of escape on a tree $T$ with uniformly bounded vertex
degrees we have
\[
{\bf HD}(\nu)=\frac{h}{c}
\]
where $c$ is the rate of escape and $h$ is the asymptotic entropy of
$P$.

This result is relevant in our context. Indeed, let $A$ be a free
basis of $F$ and let $\phi: F\to Isom(X)$ be a free, discrete and
minimal isometric action of $F$ on an $\mathbb R$-tree $X$. Then $X/F$
is a finite metric graph and $X$ is the universal cover of this graph.
Let $\Gamma(F,A)$ denote the Cayley graph of $F$ with respect to $A$.
The orbit map $w\mapsto wp$ (where $p\in X$ is a base-point) gives a
quasi-isometry between the trees $\Gamma(F,A)$ and $X$ which extends
to a homeomorphism $\hat\phi:\partial \Gamma(F,A)\to \partial X$ where
$\partial X$ is metrized in the standard $CAT(-1)$ way:
$d(\zeta,\xi)=e^{-d(p, [\zeta,\xi])}$ for $\zeta,\xi\in \partial X$.
Let $\mu_A$ be the uniform probability measure on $\partial
\Gamma(F,A)=\partial F$ corresponding to $A$ and let $\mu_A'$ denote
the push-forward of $\mu_A$ via $\hat\phi$ to $\partial X$.

Then the result of Kaimanovich~\cite{Kaim98} mentioned above implies
that
\[
{\bf HD}_d(\mu_A')=\frac{\log (2k-1)}{\lambda_A(\phi)},
\]
where $k\ge 2$ is the rank of $F$.

\subsection{Main results about actions on trees}$ $

Our first main result is:
\begin{theor}\label{free}
  Let $F=F(A)$ be a free group of rank $k\ge 2$.  Let $\phi:F\to
  Aut(X)$ be a free simplicial action without inversion of $F$ on a
  simplicial tree $X$.

  Then the following hold:

\begin{enumerate}

\item The generic stretching factor $\lambda=\lambda(\phi)$ is a
  rational number $\ge 1$ with
\[
2k\lambda\in {\mathbb Z}\left[\frac{1}{2k-1}\right].
\]

%\item The number $\lambda$ has the following additional properties:

%  (a) For any $\epsilon>0$ and for any point $p\in X$ if $P_n$ is the uniform probability
%measure on the set of elements of $F$ of length $n$ then

%\[
%\lim_{n\to\infty} P_n\left(\frac{|f|_p}{|f|}\in (\lambda-\epsilon,
%\lambda+\epsilon)\right)=1 ,
%\]
%and the convergence is exponentially fast.

%(b) For any $\epsilon>0$ if $P_n'$ is the uniform probability
%measure on the set of cyclically reduced elements of $F$ of length
%$n$ then

%\[
%\lim_{n\to\infty} P_n'\left(\frac{||\phi(w)||_p}{||w||}\in (\lambda-\epsilon,
%\lambda+\epsilon)\right)=1 ,
%\]
%and the convergence is exponentially fast.

\item The number $\lambda(\phi)$ is algorithmically computable in
  terms of $\phi$, provided $X$ is the universal cover of a finite
  connected graph and $\phi$ is given by an isomorphism between $F$
  and the fundamental group of that graph.

\end{enumerate}

\end{theor}

The most interesting case of the above theorem is when $X$ is the
Cayley graph of $F=F(A)$ determined by an endomorphism of $F$.
\begin{defn}[Generic stretching factor of an endomorphism]
  Let $F=F(A)$ where $k\ge 2$ and $A=\{a_1,\dots, a_k\}$. Let
  $\phi:F\to F$ be an endomorphism of $F$. Let $X=\Gamma(F,A)$ be the
  Cayley graph of $F$ and consider the action $\theta:F\to Isom(X)$
  given by $\theta(w)x:=\phi(w)x$, where $w\in F,x\in X$.  The generic
  stretching factor $\lambda_A(\theta)$ corresponding to this action
  is called the \emph{generic stretching factor of $\phi$ with respect
    to $A$} and is denoted $\lambda_A(\phi)$ or just $\lambda(\phi)$
  if $A$ is fixed.
\end{defn}

Thus $\lambda(\phi)$ approximates the distortion $|\phi(w)|_A/|w|_A$
for a long random freely reduced word $w$ in $A^{\pm 1}$. For
instance, for the Nielsen automorphism $\phi\in Aut(F(a,b))$,
$\phi(a)=ab, \phi(b)=b$ it turns out that $\lambda(\phi)=\frac{7}{6}$.
If $\phi$ is an automorphism of $F(a_1,\dots, a_k)$, then the precise
relationship between $\lambda(\phi)$ and the traditionally studied
dynamical properties of $\phi$ is not very clear. Nevertheless, we are
able to estimate the growth of $\lambda(\phi^n)$ for hyperbolic
automorphisms.  Recall that $\phi\in Aut(F)$ is \emph{hyperbolic} if
there exist $s >1$ and $m\ge 1$ such that for any $w\in F$
\[
s ||w||\le \max\{||\phi^m(w)||, ||\phi^{-m}(w)||\}.
\]
By a result of Brinkmann~\cite{Br} an automorphism $\phi\in Aut(F)$ is
hyperbolic if and only if $\phi$ does not have any nontrivial periodic
conjugacy classes in $F$. We prove:

\begin{theor}\label{hyp}
  Let $F=F(a_1,\dots, a_k)$ and let $\phi\in Aut(F)$ be a hyperbolic
  automorphism with parameters $s >1$ and $m\ge 1$ as above. Then
\[
\liminf_{n\to\infty}\sqrt[n]{\lambda(\phi^n)} \ge s^{1/m}>1.
\]
\end{theor}

It is obvious that any automorphism of a finitely generated group $G$
equipped with a word metric, is a quasi-isometry and indeed a
bi-Lipschitz equivalence. However, from the geometric point of view,
especially in light of various versions of the Marked Length Spectrum
Rigidity Conjecture, it is natural to study finer features of
quasi-isometries. Recall that a map $f:(X,d)\to (X',d')$ is called a
\emph{rough isometry} if there is $D>0$ such that for any $x,y\in X$
we have $|d'(f(x),f(y))-d(x,y)|\le D$. A map $f:(X,d)\to (X',d')$ is
called a \emph{rough similarity} if there are $\lambda >0$ and $D>0$
such that for any $x,y\in X$ we have $|d'(f(x),f(y))-\lambda
d(x,y)|\le D$.  It is interesting and natural to ask when an
automorphism is a rough similarity or a rough isometry.

An automorphism $\phi$ of $F=F(A)$ is called a \emph{relabelling
  automorphism} if it is induced by a permutation of the set $A =
\{a_1,\dots, a_k\}^{\pm 1}$.  We say that $\phi \in Aut(F)$ is
\emph{simple} if it is equivalent to a relabeling automorphism in
$Out(F)$, that is, if $\phi$ is the composition of a relabeling
automorphism and a conjugation. Note that being a simple automorphism
has a nice geometric meaning. Let $F=F(a_1,\dots, a_k)$ be realized as
the fundamental group of the metric graph $\Gamma$ which is a bouquet
of $k$ circles of length $1$ corresponding to the generators
$a_1,\dots, a_k$. An automorphism $\phi$ is simple if and only if,
after possibly a composition with an inner automorphism , $\phi$ is
induced by an \emph{isometry} of the graph $\Gamma$.

Let $P_n$ be the uniform probability measure on the set of all
elements of $F$ of length $n$. A set $W\subseteq F$ is said to be
\emph{exponentially $F$-generic} if $\lim_{n\to\infty} P_n(W)=1$ and
convergence to this limit is exponentially fast.  Similarly, a subset
$C \subseteq \mathcal{CR}$ of the set $\mathcal{CR}$ of all cyclically
reduced words is \emph{exponentially $\mathcal{CR}$-generic} if
$\lim_{n\to\infty} P_n'(C)=1$ with exponentially fast convergence,
where $P_n'$ is the uniform discrete probability measure on the set of
cyclically reduced words of length $n$.

Obviously, any simple automorphism is a rough isometry and a rough
similarity. The converse is also true, that is, any automorphism which
is a rough similarity must be simple. (This follows, for example, from
Theorem~2 of \cite{Fur} together with some standard results about
Culler-Vogtmann outer space). Here we obtain a strengthened "random
rigidity" version of this fact:

\begin{theor}\label{rigid}
  Let $F=F(a_1,\dots, a_k)$ be a free group of rank $k\ge 2$ with the
  standard word metric $d$ corresponding to the free basis
  $\{a_1,\dots, a_k\}$.  Put $d_0:=1+\frac{2k-3}{4k^2-2k}$. There
  exists an exponentially $\mathcal{CR}$-generic set $C \subseteq
  \mathcal{CR}$ with the following property.

  For any $\phi\in Aut(F)$ the following conditions are equivalent:

\begin{enumerate}

\item The automorphism $\phi$ is simple.

\item We have $\lambda(\phi)=1$.

\item We have $\lambda(\phi)< 1+\frac{2k-3}{2k^2-k}$.

\item The map $\phi: (F,d)\to (F,d)$ is a rough isometry.

\item The map $\phi: (F,d)\to (F,d)$ is a rough similarity.

\item For some $w\in C$ we have $||\phi(w)||=||w||$.

\item For every $w\in C$ we have $||\phi(w)||=||w||$.

\item For some $w\in C$ we have $||\phi(w)||\le d_0 ||w||$.

\item For every $w\in C$ we have $||\phi(w)||\le d_0 ||w||$.

\end{enumerate}
\end{theor}

This result shows, in particular, that the set of all possible values
of $\lambda(\phi)$ (where $\phi\in Aut(F)$) has a \emph{gap}, namely
the interval $(1, 1+\frac{2k-3}{2k^2-k})$.  Moreover, in the above
theorem we can choose $d_0$ to be any number such that
$1<d_0<1+\frac{2k-3}{2k^2-k}$.

Theorem~\ref{rigid} introduces a new dimension for rigidity results
related to Marked Length Spectra on hyperbolic groups.  Indeed, it is
well-known that if $\phi\in Aut(F)$ fixes the lengths of all conjugacy
classes (that is of all cyclic words), then $\phi$ is a rough isometry
of $F$.  Theorem~\ref{rigid} shows that even if $\phi\in Aut(F)$
preserves the length of a single ``random'' cyclically reduced word
$w$ then $\phi$ is a rough isometry and indeed a simple automorphism.
To prove Theorem~\ref{rigid} we need some rather different tools and
ideas, both algebraic and probabilistic.  The key ingredient there is
the work of Kapovich-Schupp-Shpilrain~\cite{KSS} about the behavior of
Whitehead's algorithm and the action of $Aut(F)$ on ``random''
elements of $F$.

Using Theorem~\ref{rigid} it is not hard to show that the set of
generic stretching factors taken over \emph{all} free actions of
$F(a_1,\dots, a_k)$ on simplicial trees also has a gap. Thus we
obtain:

\begin{theor}\label{rigid1}
  Let $F=F(a_1,\dots, a_k)$ where $k\ge 2$. Let $\phi: F\to Aut(X)$ be
  a free minimal action on $F$ on a simplicial tree $X$ without
  inversions.

  Then exactly one of the following occurs:

\begin{enumerate}
\item There is a simple automorphism $\alpha$ of $F$ such that $X$ is
  $\phi\circ \alpha$-equivariantly isomorphic to the Cayley graph of
  $F$ with respect to $\{a_1,\dots, a_k\}$. In this case
  $\lambda(\phi)=1$.
\item We have $\lambda(\phi)\ge 1+\frac{1}{k(2k-1)}$.
\end{enumerate}
\end{theor}

For an automorphism $\phi\in Aut(F)$ the \emph{conjugacy distortion
  spectrum} of $\phi$ is
\[
I(\phi):=\left\{ \frac{||\phi(w)||}{||w||}: w\in F-\{1\} \right\}.
\]
Kapovich proved in~\cite{Ka} that $I(\phi)$ is always a $\mathbb
Q$-convex subinterval of $\mathbb Q$ (that is, a set closed under
taking rational convex combinations) with rational endpoints.  Here we
obtain:

\begin{corol} Let $\phi\in Aut(F)$ be an arbitrary automorphism.
  Then the following hold.

\begin{enumerate}
\item Either $\phi$ is simple and $I(\phi)=1$ or, else, $1$ belongs to
  the interior of $\overline{I(\phi)}$.
\item There exists $w\in F, w\ne 1$ such that $||\phi(w)||=||w||$.
\end{enumerate}
\end{corol}
\begin{proof}
  Part (1) obviously implies part (2) since, by the above mentioned
  result of~\cite{Ka} $I(\phi)$ is a $\mathbb Q$-convex subset of
  $\mathbb Q$.

  To see that (2) holds, assume that $\phi$ is not simple. Hence
  $\phi^{-1}$ is not simple either. By Theorem~\ref{rigid} we have
  $\lambda(\phi)>1$ and $\lambda(\phi^{-1})>1$. Part (2b) of
  Theorem~\ref{free} now implies that there exists $w_1, w_2\in F,
  w_1\ne 1, w_2\ne 1$ such that $x:=\frac{||\phi(w_1)||}{||w_1||}>1$
  and $y:=\frac{||\phi^{-1}(w_2)||}{||w_2||}>1$. By definition of
  $I(\phi)$ we have $x\in I(\phi)$. Also, with $u=\phi^{-1}(w_2)$ we
  have $y=\frac{||u||}{||\phi(u)||}>1$ and so $1/y\in I(\phi)$.  Since
  $x>1$ and $1/y<1$, the $\mathbb Q$-convexity of $I(\phi)$ implies
  that $1$ belongs to the interior of $\overline{I(\phi)}$, as
  claimed.
\end{proof}

We also obtain an application of Theorem~\ref{rigid} concerning the
notion of the \emph{flux} of an automorphism that was introduced and
studied by Myasnikov and Shpilrain in~\cite{MS1}.

\begin{defn}[Flux]\cite{MS1}
  Let $G$ be a finitely generated group with a fixed word metric.  Let
  $\phi\in Aut(G)$.

  For each $n\ge 0$ define
\[
flux_{\phi}(n):=\# \{ g\in G : |g|\le n, |\phi(g)|>n\}
\]

and
\[
flux(\phi):=\limsup \sqrt[n]{\frac{flux_{\phi}(n)}{\# \{g\in G: |g|\le
    n\}}}
\]
\end{defn}

The sequence $flux_{\phi}(n)$ and the number $flux(\phi)$ provide a
certain dynamical "measure of activity" of an automorphism $\phi$. As
a corollary of our results in this paper we obtain:
\begin{corol}\label{cor:flux}
  Let $F=F(a_1,\dots, a_k)$ be a free group of rank $k\ge 2$, equipped
  with the standard metric.

  Then for any $\phi\in Aut(F)$ we have:

\[
flux(\phi)=\begin{cases} 0, \text{ if } \phi \text{ is a
    relabeling automorphism}\\
  1, \text{ otherwise.}
\end{cases}
\]
\end{corol}

\subsection{Random elements in regular languages}
By definition, a language $L$ over the alphabet $A $ is regular if and
only if there is a deterministic finite automaton which accepts the
language $L$.  It is a basic fact of formal language theory that the
class of languages accepted by nondeterministic finite automata (NDFA)
is also the class of regular languages.  (See Hopcroft and Ullman
~\cite{HU}.) Nondeterministic automata are very useful because a NDFA
accepting a language $L$ may be much smaller than any deterministic
automaton accepting $L$.  Such an automaton is not unique and choosing
some finite automaton accepting $L$ is like choosing a presentation
for a group. One can choose a ``random'' element in the regular
language $L$ is via a random walk in the transition graph of any
``suitable'' finite state automaton $M$ accepting the language $L$. We
make this precise in Section~\ref{sect:lang}, where we associate to
$M$ a finite state Markov process $M'$ with the set of states being
the set of directed edges in the transition graph $\Gamma(M)$ of $M$.
The sample space $\Omega$ of $M'$ is the set of semi-infinite
edge-paths in $\Gamma(M)$. Each path in $\Gamma(M)$ (finite or
infinite) has a label that is a word (finite or infinite) over $A$.
If $\omega\in \Omega$ is such an infinite path, we denote by
$w_n=w_n(\omega)$ the label of the initial segment of length $n$ of
$\omega$. Any initial probability distribution $\mu$ on the edge-set
$E(\Gamma(M))$ defines a probability measure $P_{\mu}$ on $\Omega$. We
need to impose a natural assumption on $M$ in order to guarantee that
the Markov process $M'$ is irreducible. This technical assumption,
which is frequently satisfied in practice, is made precise in the
definition of a \emph{normal} automaton in Section~\ref{sect:lang}.
Again applying the Subadditive Ergodic Theorem, we have:

\begin{theor}\label{lang}
  Let $M$ be a normal automaton over a finite alphabet $A$ and let
  $L=L(M)$ be the language accepted by $M$.

  Let $\phi: A^*\to G$ be a monoid homomorphism, where $G$ is a group
  with a left-invariant semi-metric $d_G$. Then there exists a number
  $\lambda=\lambda (M,\phi, d_G)\ge 0$ such that for any initial
  distribution $\mu$ on $E(\Gamma(M))$ we have
\[
\frac{|\phi(w_n)|_G}{n} \to \lambda \text{ almost surely and in } L^1
\text{ with respect to } P_{\mu}.
\]
\end{theor}
If the initial distribution $\mu$ is supported on the edges of
$\Gamma(M)$ originating at the start states of $M$ then the word $w_n$
can be extended by a word of uniformly bounded length to get a word
$w_n'\in L$. We can think of $w_n'$ as a "random" element of $L$ with
respect to $M$ and $\mu$. Theorem~\ref{lang} then implies that
$\frac{|\phi(w_n')|_G}{n} \to \lambda$ as $n\to\infty$ almost surely
and in $L^1$ with respect to $\mu$.

\section{Random words and random walks}\label{sect:erg}

\begin{conv}
  Let $A=\{a_1,\dots, a_k\}$ be a \emph{free generating set} of a
  \emph{free group} $F=F(A)$ of finite rank $k>1$. For $w\in F$ we
  denote by $|w|_A$ or simply $|w|$) the \emph{freely reduced length}
  of $w$ with respect to $A$. Let $d(w_1,w_2)=\left| w_1^{-1} w_2
  \right|$ be the associated left-invariant metric on $F$. By
  $||w||_A=||w||$ we denote the \emph{cyclically reduced length} of
  $f$ with respect to $A$, that is, the length of any cyclically
  reduced word in the alphabet $A^{\pm 1}$ conjugate to $f$.

%The \emph{hyperbolic boundary} $\partial F$ (which, in this case, coincides
%with the \emph{space of ends} of $F$) is naturally identified with the set of
%\emph{semi-infinite freely reduced words} in $A^{\pm 1}$ corresponding to
%geodesic rays from the identity $1$ in the standard Cayley graph $\Gamma(F,A)$.

%Denote by $\omega_n\in F$ the point on the ray $\omega\in\partial
%F$ at distance $n$ from the group identity, and for a freely
%reduced word $u$ of length $n$ let
%$$
%Cyl_A(u) = \{\omega\in\partial F: \omega_n =u \}
%$$
%be the set of all $\omega\in\partial F$ that have $u$ as initial segment. The
%\emph{cylinder sets} $Cyl_A(u)$ generate the Borel $\sigma$-algebra $\mathcal F$
%in $\partial F$.

%The \emph{uniform} Borel probability measure $\m=\m(A)$ on $\partial F$ is
%defined by assigning equal weights to all cylinders based on the words $u$ of
%the same length $n=|u|$, i.e.,
%\[
%\m(Cyl_A(u))=\begin{cases} 1  &\;, \text{ if } n=0 \\
%  {\displaystyle\frac{1}{2k(2k-1)^{n-1}}} &\;, \text{ if } n>0
% \end{cases}
%\]

  This convention, including the fixed choice of the free basis
  $A=\{a_1,\dots, a_k\}$ of $F$, is adopted for the remainder of the
  paper, unless specified otherwise.
\end{conv}

Recall that a nonnegative function $|\cdot|_G$ on a group $G$ is
called a \emph{semi-norm} if for all $g,h\in G$ we have $|gh|_G\le
|g|_G+|h|_G$.

In this Section we shall prove Theorem~\ref{thm:stretch} from the
Introduction:

\begin{thm} \label{thm:rate of escape}
  Let $F=F(A)$, and let $\m=\m(A)$ be the uniform Borel probability
  measure on $\partial F$ corresponding to the generating set
  $A=\{a_1,\dots, a_k\}$ with $k\ge 2$. Let $\phi:F\to G$ be a
  homomorphism to a group $G$ endowed with a semi-norm $|\cdot|_G$.
  Then:

\begin{enumerate}
\item There exists a real number $\lambda\ge 0$ such that
\[
\lim_{n\to\infty} \frac{|\phi(\omega_n)|_G}{n}=\lambda
\]
for $\m$-a.e. $\omega\in\partial F$ and in the space $L^1(\partial
F,\m)$.

\item If the group $\phi(F)$ is non-amenable, the sequence
\[b_n=\#\{ g\in \phi(F): |g|_G\le n\}\] grows at most exponentially,
then $\lambda>0$.
\end{enumerate}
\end{thm}
The condition on $b_n$ in the above theorem is always satisfied if $G$
is a finitely generated group and $|.|_G$ is the word metric
corresponding to some finite generating set of $G$.

Any geodesic ray $\omega\in\partial F$ can be identified with the
\emph{non-backtracking path} $\omega_0,\omega_1,\dots$ in $F$ starting
from the group identity. Then the measure space $(\partial F,\m)$
becomes the \emph{space of sample paths} of the \emph{non-backtracking
  simple random walk} (NBSRW) on the Cayley graph of $F$ starting from
the identity of the group.  This is the Markov chain on $F$ whose
transition probabilities $\pi_f,\,f\in F$ are equidistributed among
the neighbors of $f$ which are strictly further from the group
identity. Therefore, the number $\lambda$ above is the \emph{linear
  rate of escape} of the $\phi$-image of the non-backtracking simple
random walk on $F$. We shall deduce Theorem~\ref{thm:rate of escape}
from well-known analogous properties of the usual random walks on
groups by using the fact that the simple random walk on the free group
asymptotically follows uniformly distributed geodesic rays.

\medskip

Let $\mu$ be a probability measure on a group $G$. By definition, the
sample paths of the associated \emph{random walk} $(G,\mu)$ are
products $g_n=h_1 h_2 \dots h_n$ of independent $\mu$-distributed
\emph{increments} $h_n$. In other words, the measure $\mathbf P$ in
the space of sample paths which describes the random walk $(G,\mu)$ is
the image of the product measure $\mu\otimes\mu\otimes\dots$ in the
space of increments under the above product map.

The following statement is known as Kingman's Subadditive Ergodic
Theorem~\cite{King}. (See also~\cite{Der} for a short proof.)

\begin{prop}[Subadditive Ergodic Theorem]\label{thm:set}
  Let $(\Omega, {\mathcal F}, \mu)$ be a probability space and let
  $\mathcal S:\Omega\to\Omega$ be a measure-preserving operator, that
  is such that for any measurable set $Q\subseteq \Omega$ we have
  $\mu(Q)=\mu(\mathcal S^{-1}Q)$.

  Let $X_n:\Omega\to \mathbb R$ be a sequence of non-negative
  integrable random variables such that for any $n,m\ge 0$
\[
X_{n+m}(\omega) \le X_n(\omega)+ X_m (\mathcal S^n\omega), \quad
\text{ a. e. } \omega\in \Omega.
\]
Then there exists a $\mathcal S$-invariant random variable
$\lambda:\Omega\to\mathbb R$ such that \[ \lim_{n\to\infty}
\frac{X_n}{n} =\lambda
\]
almost surely and in $L^1$ on $\Omega$.

In particular, if $\mathcal S$ is ergodic then $\lambda=const$ on
$\Omega$.
\end{prop}

A straightforward application of Kingman's Subadditive Ergodic Theorem
gives:

\begin{prop}[\cite{G}] \label{prop:rate}
  If the measure $\mu$ has a finite first moment $\sum |g| \mu(g)$
  with respect to a semi-norm $|\cdot|$ on the group $G$, then there
  exists a number $c\ge 0$ (called the linear rate of escape of the
  random walk $(G,\mu)$ with respect to the semi-norm $|\cdot|$) such
  that $|g_n|/n\to c$ for $\mathbf P$-a.e.  sample path $(g_n)$ and in
  the space $L^1(\mathbf P)$.
\end{prop}

The following claim, if slightly more general than the one formulated
in \cite{G}, can be obtained in the same way by using the spectral
characterization of amenability (or by showing that $c=0$ implies
vanishing of the asymptotic entropy of the random walk, and therefore
amenability of the group \cite{KV}):

\begin{prop}[\cite{G}] \label{prop:positive rate}
  Under the assumptions of Proposition~\ref{prop:rate}, if the group
  $G$ is non-amenable and the semi-norm $|\cdot|_G$ has exponentially
  bounded growth and the support of the measure $\mu$ generates the
  group $G$, then $c >0$.
\end{prop}

Let now $\mu_A'$ be the probability measure on the free group $F$
equidistributed on the set $A^{\pm 1}$, so that $\mu_A'\left(a_i^{\pm
    1}\right)=1/2k$ for $i=1,2,\dots,k$.

\begin{prop}[see \cite{K} and the references therein] \label{prop:asymptotic}
  For $\mathbf P$-a.e. sample path $(g_n)$ of the random walk
  $(F,\mu_A')$

\begin{enumerate}

\item There exists a limit
\[
g_\infty=\lim_n g_n\in\partial F,
\]
and its distribution (i.e., the image of the measure $\mathbf P$ under
the map $(g_n)\mapsto g_\infty$) coincides with the uniform measure
$\m$ on $\partial F$.

\item We have
\[
\lim_n \frac{|g_n|}{n} = \theta = \frac{k-1}{k}.
\]
(so that the linear rate of escape of the random walk $(F,\mu_A')$ is
$\frac{k-1}{k}$).

\item We have
\[
d \bigl( g_n, (g_\infty)_{[\theta n]} \bigr) = o(n),
\]
where $[x]$ denotes the integer part of a number $x$.
\end{enumerate}

\end{prop}

\begin{proof}[Proof of Theorem~\ref{thm:rate of escape}]
  Consider the random walk $(F,\mu_A')$. Its image under the
  homomorphism $\phi$ is the random walk on the group $\phi(F)$
  governed by the measure $\phi(\mu_A')$. Denote its rate of escape
  with respect to the semi-norm $|\cdot|_G$ by $c$. Then the
  combination of Proposition~\ref{prop:rate} and
  Proposition~\ref{prop:asymptotic} implies the first part of
  Theorem~\ref{thm:rate of escape}. Indeed, for $\mathbf P$-a.e.
  sample path $(g_n)$ the distance in $F$ between $g_n$ and
  $(g_\infty)_{[\theta n]}$ is sublinear, whence the distance in $G$
  (with respect to the semi-metric determined by the semi-norm
  $|\cdot|_G$) between the $\phi$-images of these points is also
  sublinear. Since the distribution of $g_\infty$ is $\m$, we arrive
  at the conclusion that the first part of Theorem~\ref{thm:rate of
    escape} holds for the number $\lambda=c/\theta$.

  The second part of Theorem~\ref{thm:rate of escape} is now an
  immediate corollary of Proposition~\ref{prop:positive rate}.

  Here is another argument establishing part (1) of
  Theorem~\ref{thm:rate of escape} as a direct consequence of the
  Subadditive Ergodic Theorem.

  Let $\Omega=\partial F$. Recall that for $\omega\in \partial F$ we
  denote by $\omega_n$ the element of $F$ that is at distance $n$ from
  $1$ along the geodesic ray $\omega$ in $\Gamma(F,A)$.  Let
  $X_n:\partial F\to\mathbb R$ be defined as
  $X_n(\omega):=|\phi(\omega_n)|_G$.  Also, let $\mathcal S:\partial
  F\to \partial F$ be the standard shift operator consisting in
  erasing the first letter of a semi-infinite freely reduced word
  representing an element of $\partial F$.  It is well-known that
  $\mathcal S$ is stationary and ergodic.

  Note that for any $\omega\in \partial F$ we have
\[
\omega_{n+m}=\omega_n (\mathcal S^n\omega)_m.
\]
Hence
\[
|\phi(\omega_{n+m})|_G= |\phi(\omega_n) \phi((\mathcal
S^n\omega)_m)|_G\le |\phi(\omega_n)|_G+ |\phi((\mathcal
S^n\omega)_m)|_G.
\]
Thus the conditions of the Subadditive Ergodic Theorem are satisfied
and part (1) of Theorem~\ref{thm:rate of escape} follows.
\end{proof}

The following is Theorem~\ref{aut} from the Introduction.

\begin{thm}\label{auto}
  Let $G$ be a nontrivial finitely generated group with a word-metric
  $d_S$ corresponding to a finite generating set $S$.  Let $H \le
  Aut(G)$ be a noncyclic finitely generated group with a finite
  generating set $T$. Then:

\begin{enumerate}
\item There is $\nu=\nu(H)=\nu(H,T,S)\ge 0$ such that for a
  non-backtracking simple random walk $\phi_n$ on the Cayley graph of
  $H$ with respect to $T$ we have

\[
\lim_{n\to\infty}\frac{\log ||\phi_n||_S}{n}=\nu
\]
almost surely and in $L^1$.

\item If $G$ has polynomial growth and $H$ is non-amenable then
  $\nu(H,T,S)>0$.
\end{enumerate}
\end{thm}

\begin{proof}
  It is clear from the definition of $||\cdot||_S$ that for any
  $\phi,\psi\in Aut(G)$ we have
\[
||\phi\psi||_S\le ||\phi||_S ||\psi||_S
\]
and hence
\[
\log ||\phi\psi||_S\le \log ||\phi||_S +\log ||\psi||_S.
\]
Also, for any $\phi\in Aut(G)$ we have $||\phi||_S\ge 1$ and so $\log
||\phi||_S\ge 0$. Thus $\log ||\cdot||_S$ is a semi-norm on $Aut(G)$
that uniquely extends to a left-invariant semi-metric of $Aut(G)$ and
thus on $H\le Aut(G)$. Hence part (1) of Theorem~\ref{auto} follows
directly from part (1) of Theorem~\ref{thm:stretch}.

To see that part (2) holds suppose that $H$ is nonamenable and that
$G$ has polynomial growth. This implies that $(H, ||\cdot||_S)$ has at
most exponential growth. Hence part (2) of Theorem~\ref{auto} follows
from part (2) of Theorem~\ref{thm:stretch}.
\end{proof}

\begin{rem}
  The requirement of $G$ having polynomial growth in Theorem~\ref{aut}
  is important and cannot be easily dispensed with. If $G$ is a
  group and $g\in G$, denote by $ad(g)\in Aut(G)$ the inner
  automorphism of $G$ defined as $ad(g)(x)=gxg^{-1}$ for every
  $x\in G$. Now let $G=F(a_1,\dots, a_k)$ and $H= Inn(F)\le Aut(F)$ be the
  (non-amenable!) group of inner automorphisms of $F$ with the
  generating set $T=\{ad(a_1), \dots, ad(a_k)\}$.
  Then for any product $\phi_n$ of $n$ elements
  of $T$ we have $||\phi_n||\le 2n+1$. Since $\lim_{n\to\infty}
  \frac{\log 2n+1}{n}=0$, we see that $\nu(H,T,S)=0$. Nevertheless, in
  some instances quotient group considerations still imply that
  $\nu(A)>0$ even if $G$ does not have polynomial growth, or,
  equivalently, $G$ is not virtually nilpotent.
\end{rem}

We obtain Corollary~\ref{cor:aut} from the Introduction:

\begin{cor}\label{cor:auto}
  Let $F$ be a free group of finite rank $k>1$ and let $H\le Aut(F)$
  be a finitely generated group of automorphisms of $F$ such that the
  image $H'$ of $H$ in $Aut(F_{ab})\cong GL(k,{\mathbb Z})$ is
  nonamenable. Then for any finite generating set $S$ of $F$ and for
  any finite generating set $T$ of $H$ we have $\nu(H,T,S)>0$.
\end{cor}

\begin{proof}
  Let $S'$ be the image of $S$ in the abelianization ${\mathbb
    Z}^k=F_{ab}$ of $F$. For any $\phi\in Aut(F)$ the automorphism
  $\phi$ of $F$ factors through to an automorphism $\phi'$ of
  $F_{ab}$. Clearly $||\phi||_S\ge ||\phi'||_{S'}$. Hence
  $\nu(H,T,S)\ge \nu(H',T',S')$, where $T'$ is the image of $T$ in
  $Aut(F_{ab})$. Since $F_{ab}$ has polynomial growth, by
  Theorem~\ref{aut} we have $\nu(H',T',S')>0$ and hence
  $\nu(H,T,S)>0$.
\end{proof}

\section{Frequencies and cyclic words}\label{sect:freq}

The following convention is fixed until the end of the paper unless
specified otherwise.

\begin{conv}

  As before, let $k\ge 2$ and let $F=F(A)$ where $A=\{a_1,\dots,
  a_k\}$. Let $\widehat{A}=A^{\pm 1}$. We denote by $\mathcal{CR}$ the
  set of all cyclically reduced words in $F$.

  A \emph{cyclic word} is an equivalence class of nontrivial
  cyclically reduced words, where two nontrivial cyclically reduced
  words are equivalent if they are cyclic permutations of each other.
  If $v$ is a cyclically reduced word, we denote by $(v)$ the cyclic
  word defined by $v$. Recall that if $u$ is a freely reduced word, we
  denote the length of $u$ by $|u|$ and the length of the cyclically
  reduced form of $u$ by $||u||$. If $w=(v)$ is a cyclic word then
  $||w||=||v||$ is the length of $w$.
\end{conv}

Note that the set of cyclic words is naturally identified with the set
of nontrivial conjugacy classes of elements of $F$.

\begin{defn}[Frequencies]
  Let $w$ be a cyclic word.

  Let $v$ be a nontrivial freely reduced word.  We define $n_w(v)$,
  \emph{the number of occurrences of $v$ in $w$} as follows.  Let
  $w=(z)$. Take the smallest $p>0$ such that $|z^{p-1}|\ge 2|v|$ and
  count the number of those $i, 0\le i<||w||$ such that $z^p\equiv
  z_1vz_2$ where $|z_1|=i$. This number by definition is $n_w(v)$. If
  $v=1$, we define $n_w(1):=||w||$.

  There is a more graphical way of defining $n_w(v)$ for a nontrivial
  cyclic word $w$.  We will think of $w$ as a cyclically reduced word
  written on a circle in a clockwise direction without specifying a
  base-point.  Then $n_w(v)$ is the number of positions on the circle,
  starting from which it is possible to read the word $v$ going
  clockwise along the circle (and wrapping around more than once, if
  necessary).

  For any freely reduced word $v$ we define \emph{frequency} of $v$ in
  $w$ as:

\[
f_w(v):=\frac{n_w(v)}{||w||}.
\]

Also, if $w$ is a nontrivial freely reduced word, and $v$ is another
nontrivial freely reduced word, we define $n_w(v)$, \emph{the number
  of occurrences of $v$ in $w$}, as follows. If $|w|=n>0$ then by
definition $n_w(v)$ is the number of those $i, 0\le i<n$ for which $w$
decomposes as a freely reduced product $w=w'vw''$ with $|w'|=i$. Thus
if $|v|\le |w|$ then necessarily $n_w(v)=0$ (unlike the situation when
$w$ is a cyclic word).
\end{defn}

\begin{lem}\label{red}
  Let $w$ be a nontrivial cyclic word. Then:

\begin{enumerate}

\item For any $m\ge 0$ and for any freely reduced word $u$ with
  $|u|=m$ we have:
\[
n_w(u)=\sum_{x\in \widehat{A},|ux|=|u|+1} n_w(ux)=\sum_{x\in
  \widehat{A},|xu|=|u|+1} n_w(xu),
\]
and
\[
f_w(u)=\sum_{x\in \widehat{A},|ux|=|u|+1} f_w(ux) =\sum_{x\in
  \widehat{A},|xu|=|u|+1} f_w(xu).
\]

\item For any $m\ge 1$
\[
\sum_{|u|=m} n_w(u)=||w|| \quad \text{ and } \sum_{|u|=m} f_w(u)=1.
\]

\item For any $s>0$ and any $u\in F$
\[
n_{w^s}(u)=sn_w(u) \quad \text{ and }\quad f_{w^s}(u)=f_w(u).
\]
\end{enumerate}

\end{lem}

\begin{proof}
  Parts (1) and (3) are obvious. We establish (2) by induction on $m$.
  For $m=1$ the statement is clear. Suppose that $m>1$ and that (2)
  has been established for $m-1$.

  We have:

\[
\sum_{|u|=m} n_w(u) = \sum_{\overset{|v|=m-1, x\in \widehat{A}:}
  {|vx|=m}} n_w(vx) = \sum_{|v|=m-1} n_w(v)=||w||,
\]
as required.

\end{proof}

\begin{defn}[Nielsen automorphisms]
  A \emph{Nielsen} automorphism of $F$ is an automorphism $\tau$ of
  one of the following types:

\begin{enumerate}

\item There is some $i, 1\le i\le k$ such that $\tau(a_i)=a_i^{-1}$
  and $\tau(a_j)=a_j$ for all $j\ne i$.

\item There are some $1\le i<j\le k$ such that $\tau(a_i)=a_j$,
  $\tau(a_j)=a_i$ and $\tau(a_l)=a_l$ when $l\ne i, l\ne j$.

\item There are some $1\le i<j\le k$ such that $\tau(a_i)=a_ia_j$ and
  $\tau(a_l)=a_l$ for $l\ne i$.

\end{enumerate}

\end{defn}

It is a classical fact that the set of all Nielsen automorphisms
generates $Aut(F)$.

The following proposition proved by Kapovich in \cite{Ka} is crucial
for our arguments.

\begin{prop}\label{phi}
  Let $\phi\in Out(F)$ be an outer automorphism and let $p\ge 0$ be
  such that $\phi$ can be represented, modulo $Inn(F)$, as a product
  of $p$ Nielsen automorphisms.

  Then for any freely reduced word $v\in F$ with $|v|=m$ there exists
  a collection of computable integers $c(u,v)=c(u,v,\phi)\ge 0$, where
  $u\in F$, $|u|=8^p m$, such that for any nontrivial cyclic word $w$
  we have

\[
n_{\phi(w)}(v)=\sum_{|u|=8^p m} c(u,v) n_w(u).
\]
\end{prop}

\begin{cor}\label{use'}
  Let $\phi$ be an automorphism of $F$ and let $p$ be such that $\phi$
  can be written as a product of $p$ Nielsen automorphisms.

  There exists a collection of integers $e(v)=e(v,\phi)\ge 0$, where
  $v\in F, |v|=8^p$, such that for any cyclic word $w$ we have:

\[
||\phi(w)||=\sum_{|v|=8^p} e(v) n_w(v)
\]
and

\[
\frac{||\phi(w)||}{||w||}=\sum_{|v|=8^p} e(v) f_w(v).
\]

Moreover, there is an algorithm which, given $\phi$ and $u$, computes
the numbers $e(v)$.
\end{cor}
\begin{proof}
  Since $\displaystyle ||\phi(w)||=\sum_{x\in \widehat{A}}
  n_{\phi(w)}(x)$, the statement follows directly from
  Proposition~\ref{phi}.
\end{proof}

The following well-known fact is a version of the so-called ``Bounded
Cancellation Lemma'' (see \cite{Coo}):

\begin{lem}\label{bnd}
  Let $\alpha$ be an injective endomorphism of $F$. There is
  $N=N(\alpha)>0$ such that for any cyclically reduced word $w$ the
  maximal terminal segment of $\alpha(w)$ that cancels in the product
  $\alpha(w)\cdot \alpha(w)$ has length at most $N$.
\end{lem}

\section{Actions on trees}

Let $\Gamma$ be a finite connected graph with orientation
$E\Gamma=E^+\Gamma\sqcup E^-\Gamma$. For $e\in E$ we use the following
notation.  The inverse edge of $e$ is denoted by $\overline e$, $o(e)$
denotes the initial vertex of $e$ and $t(e)$ denotes the terminal
vertex of $e$.

Let $F$ be a free group and let $\phi:F\to \pi_1(\Gamma,p)$ be an
isomorphism between $F$ and the fundamental group of $\Gamma$ with
respect to a vertex $p$. Let $T$ be a maximal tree in $\Gamma$.  For
any vertex $v$, let $[p,v]_T$ denote the path in $T$ from $p$ to $v$.
The choice of $T$ define a basis $S_T$ of $\pi_1(\Gamma,p)$ as
follows:

\[
S_T:=\{ [p,o(e)]_T~e~[t(e),p]_T: e\in E^+(\Gamma-T)\}
\]
The $\phi$-pullback of this basis $B_T:=\phi^{-1}(S_T)$ is a basis of
$F$ referred to as the \emph{geometric basis} of $F$ determined by
$T$.

Let $s_e:=[p,o(e)]_T~e~[t(e),p]_T$ where $e\in E (\Gamma-T)$, so that
$s_{\bar e}=s_e^{-1}$.  Let $b_e=\phi^{-1}(s_e)$, where $e\in E
(\Gamma-T)$, so that again $b_{\bar e}=b_e^{-1}$.

The following obvious lemma indicates the explicit correspondence
between freely reduced words in $S_T$ (or $B_T$) and reduced
edge-paths in $\Gamma$.

\begin{lem}\label{lem:rewr}
\begin{enumerate}
\item Let $\gamma$ be an edge-path in $\Gamma$ from $p$ to $p$. Let
  $u$ be a word in $S_T$ constructed from $\gamma$ as follows: delete
  all the edges of $T$ from $\Gamma$ and replace each edge $e\in E^+
  (\Gamma-T)$ in $\gamma$ by $s_e$ and each edge $e\in E^- (\Gamma-T)$
  in $\gamma$ by $s_{\bar e}^{-1}$.  Then $u=\gamma$ in
  $\pi_1(\Gamma,p)$ and $u$ is a reduced word in $S_T$ if and only if
  $\gamma$ is a reduced path.  x\item Let $u$ be a word in $S_T\cup
  S_T^{-1}$, where $\epsilon_i=\pm 1$.

  Construct the path $\gamma$ from $p$ to $p$ as follows. First for
  each $e\in E^+(\Gamma-T)$ replace every $s_e$ in $u$ by $e$ and
  replace every $s_e^{-1}$ by $\bar e$. Then between every two
  consecutive $e,e'$ insert the path $[t(e),o(e')]_T$. Finally append
  the path $[p,o(e)]_T$ in front, for the first edge $e_0\in
  E(\Gamma-T)$ of the resulting sequence, and append the path
  $[t(e_0'),p]_T$ at the end for the last edge $e_0'\in E(\Gamma-T)$
  of the sequence.

  Then $\gamma$ is a path from $p$ to $p$ that is equal to $u$ in
  $\pi_1(\Gamma,p)$ and that is reduced if and only if the word $u$
  over $S_T$ is reduced.
\item Let $\gamma$ be a closed edge-path in $\Gamma$. Let $u$ be a
  word in $S_T^{\pm 1}$ obtained from $\gamma$ as in (1). Then the
  loop at $p$ corresponding to $u$ in $\pi_1(\Gamma,p)$ is freely
  homotopic to $\gamma$ in $\Gamma$ and the word $u$ is cyclically
  reduced if and only if the path $\gamma$ is cyclically reduced.

\item Let $w$ be a cyclic word in $S_T^{\pm 1}$. Let $\gamma$ be a
  circuit in $\Gamma$ obtained as follows. Replace each occurrence of
  $s_e$ in $w$ by $e$ and each occurrence of $s_e^{-1}$ by $\bar e$;
  after that between each two consecutive (in the cyclic order) edges
  $e,e'$ insert the path $[t(e),o(e')]_T$.  Then $w$ and $\gamma$
  represent freely homotopic loops in $\Gamma$ and the cyclic word $w$
  is reduced if and only if the circuit $\gamma$ is reduced.

\end{enumerate}
\end{lem}

Now suppose that $\Gamma$ is endowed with the structure of a
\emph{metric graph}, that is, each edge $e$ of $\Gamma$ is assigned a
\emph{length} $\ell(e)>0$ in such a way that $\ell(e)=\ell(\bar e)$
for each edge $e$. Let $X=\widetilde{(\Gamma,p)}$ be the universal
cover of $\Gamma$. Then $X$ inherits the structure of a metric tree
with an isometric action of $\pi_1(\Gamma,p)$ and, via $\phi$, an
action of $F$ on $X$.

Let $\tilde p$ be a lift of $p$ to $X$. For $g\in \pi_1(\Gamma,p)$ let
$|g|_p:=d_X(\tilde p, g \tilde p)$. Also denote by $||g||_X$ the
translation length of $g$ when acting on $X$. Similarly, if $w$ is a
conjugacy class (or a cyclic word) in $\pi_1(\Gamma,p)$, we denote by
$||w||_X$ the translation length of $u$ with respect to $X$. For each
freely reduced word $z=s_{e}^{\epsilon} s_{e'}^{\delta}$ of length two
in $S_T^{\pm 1}$, where $\epsilon, \delta \in \{1,-1\}$, denote by
$r_z$ the length of the edge-path $[t({e}^{\epsilon}),
o({e'}^{\delta})]_T$ in $\Gamma$. Let $Z$ be the set of all freely
reduced words of length two in $S_T$. For each $a=s_{e}\in S_T$ denote
$e(a):=$ and $e(a^{-1}):=\bar e$.

\begin{lem}\label{lem:compute}
  (1) Let $w$ be a reduced cyclic word in $S_T^{\pm 1}$. Then
\[
||w||_X=\sum_{a\in S_T^{\pm 1}} \ell(e(a))n_w(a) + \sum_{z\in Z} r_z
n_w(z).
\]

(2) Let $u$ be a freely reduced word in $S_T$. Then
\[
|u|_p=\sum_{a\in S_T^{\pm 1}} \ell(e(a))n_u(a) + \sum_{z\in Z} r_z
n_u(z)+\ell([p,o(e)]_T)+\ell([t(e'),p]_T)
\]
where $e$ and $e'$ are the last and the first edges of $\gamma(u)$
accordingly.
\end{lem}

The following is Theorem~\ref{free} from the Introduction:

\begin{thm}\label{main}
  Let $F=F(a_1,\dots, a_k)$, where $k\ge 2$, and let $A=\{a_1,\dots,
  a_k\}$. Then for any free action $\phi$ of $F$ on a simplicial tree
  $X$ without inversions the generic stretching factor
  $\lambda(\phi)=\lambda_A(\phi)$ is a rational number with
\[
2k\lambda(\phi)\in \mathbb Z[\frac{1}{2k-1}].
\]

Moreover, if $X$ is given as the universal cover of a finite connected
simplicial graph $\Gamma$ and if the action $\phi$ is given via an
explicit isomorphism between $F$ and $\pi_1(\Gamma,p)$, then
$\lambda(\phi)$ is algorithmically computable in terms of $\phi$.
\end{thm}

\begin{proof}
  Recall that the definition of $\lambda(\phi, |\cdot|_x)$ does not
  depend on the choice of a point $x\in X$. Hence we may assume that
  $x$ is a vertex of the minimal $F$-invariant subtree of $X$ and
  therefore, that the action of $F$ on $X$ is minimal. Let
  $\Gamma=X/F$ be the finite quotient graph. Choose an orientation on
  $\Gamma$, a maximal tree $T$ in $\Gamma$. Choose a base-vertex $p$
  in $\Gamma$ to be the image of $x\in X$ in $\Gamma$. Note that in
  both $X$ and $\Gamma$ every edge has length $1$. Then there is a
  canonical isomorphism $\psi:F\to \pi_1(\Gamma,p)$. Let $S_T$ and
  $B_T$ be the geometric bases corresponding to $T$ for
  $\pi_1(\Gamma,p)$ and $F$ accordingly.

  Fix a bijection between $B_T$ and $A=\{a_1,\dots, a_k\}$ and an
  automorphism $\alpha$ of $F$ induced by this bijection of the two
  free bases of $F$.

  Let $g=x_1\dots x_n\in F$ be a freely reduced word of length $n$ in
  $F(a_1,\dots, a_k)$.  Let $g'$ be a cyclically reduced word of
  length $n$ over $A$ obtained from $g$ by changing the last letter of
  $g$ if necessary. Thus $|g'g^{-1}|_A\le 2$.

  Let $w'$ be the cyclic word over $A$ defined by $g'$. Let $w$ be the
  result of rewriting $w'$ as the cyclic word in $B_T$. Then there is
  an integer $M\ge 1$ such that for each freely reduced word $z$ in
  $B_T$ of length at most $2$
\[
n_{w}(z)=\sum_{|u|_A=M} c(u,z) n_{w}(u)
\]
where $c(u,z)\ge 0$ are some integers independent of $w$.  Let $Z_i$
be the set of freely reduced words of length $i$ in $B_T$, for
$i=1,2$.

Then
\begin{gather*}
  ||g'||_X=||w||_X=\sum_{a\in Z_1} \ell(e(a))n_w(a) + \sum_{z\in Z_2} r_z n_w(z)=\\
  \sum_{a\in Z_1} \sum_{|u|_A=M} \ell(a) c(u,a) n_{w'}(u) + \sum_{z\in
    Z_2} \sum_{|u|_A=M} r_z c(u,z) n_{w'}(u),
\end{gather*}

It follows from Lemma~\ref{lem:rewr} and Lemma~\ref{lem:compute} that
if $h\in F$ is cyclically reduced over $B_T$ then $\big|
||h||_X-|h|_x\big|\le N$ where $N=N(x)>0$ is some constant independent
of $h$.  On the other hand by the Bounded Cancellation Lemma (
Lemma~\ref{bnd}) there exists a constant $N'>0$ such that for any
cyclically reduced word $y$ over $A$, we have $\big|
||y||_{B_T}-|y|_{B_T}\big|\le N'$. By construction $g'$ is cyclically
reduced over $A$ and $|g'g^{-1}|_A\le 2$. Hence there exists a
constant $L>0$ such that for every $g$ as above and each $u\in F$ with
$|u|_A=M$ we have $\big| |g|_x-||g'||_X\big|\le L$ and
$|n_g(u)-n_{w'}(u)|\le L$.

Therefore there is another constant $L'>0$ independent of $f$ such
that for every freely reduced word $g$ of length $n$ over $A$

\[
\left|\sum_{a\in Z_1} \sum_{|u|_A=M} \ell(e(a)) c(u,a) f_{g}(u) +
  \sum_{z\in Z_2} \sum_{|u|_A=M} r_z c(u,z)
  f_{g}(u)-\frac{|g|_p}{n}\right|\le \frac{L'}{n} \tag{*}.
\]

If $g_n$ is a freely reduced word of length $n$ obtained by a
non-backtracking simple random walk of length $n$ on $F(a_1,\dots,
a_k)$ then for each $u\in F(a_1,\dots, a_k)$ with $|u|_A=M$ we have

\[
\lim_{n\to \infty} f_{g_n}(u)=\frac{1}{2k (2k-1)^{M-1}} \quad \text{
  almost surely.}
\]

Therefore (*) implies that

\[
\lambda(\phi)= \frac{1}{2k (2k-1)^{M-1}} \big[ \sum_{a\in Z_1}
\sum_{|u|_A=M} \ell(e(a)) c(u,a) + \sum_{z\in Z_2} \sum_{|u|_A=M} r_z
c(u,z) \big] \tag{**}
\]
Since $\ell(e(a))=1, c(u,a), c(u,z)$ and $r_z$ are integers, it
follows that $\lambda(\phi)$ is rational and, moreover, that
\[
2k \lambda(\phi) \in {\mathbb Z}[\frac{1}{2k-1}].
\]

Moreover, $\lambda(\phi)$ is computable in terms of an explicit
isomorphism between $F$ and $\pi_1(\Gamma,p)$.
\end{proof}

\begin{rem}
  The formula (**) for $\lambda(\phi)$ holds for an arbitrary
  structure of a metric graph on $\Gamma$, where the lengths of edges
  are allowed to be arbitrary positive real numbers and not
  necessarily $1$. If the length of all edges of $\Gamma$ are
  rational, then by (**) $\lambda(\phi)$ is also rational.  Moreover,
  if these length of the edges are given to us in some algorithmically
  describable form then $\lambda(\phi)$ is computable in terms of
  these lengths and of an an explicit isomorphism between $F$ and
  $\pi_1(\Gamma,p)$.
\end{rem}

\section{Genericity}\label{sect:gen}

\begin{conv}
  Recall that $\mathcal{CR}$ denotes the set of all cyclically reduced
  words in $F=F(a_1,\dots, a_k)$. If $S\subseteq F$ and $n\ge 0$ we
  denote
\[
\rho(S,n):=\#\{w\in S : |w|\le n\},
\]
and
\[
\gamma(S,n):=\#\{w\in S : |w|= n\},
\]
Let $P_n$ be the uniform discrete probability measure on the set of
all elements $w\in F$ with $|w|=n$. We extend $P_n$ to $F$ by setting
$P_n(w)=0$ for any $w\in F$ with $|w|\ne n$.

Similarly, let $P_n'$ be the uniform discrete probability measure on
the set of all cyclically reduced elements $w\in F$ with $||w||=n$. We
extend $P_n$ to $\mathcal{CR}$ by setting $P_n'(w)=0$ for any $w\in
\mathcal{CR}$ with $||w||\ne n$.
\end{conv}

Thus $\gamma(n,F)=2k(2k-1)^{n-1}$ for $n>0$.

For a number sequence $x_n$ with $\lim_{n\to\infty} x_n=x\in \mathbb
R$ we say that the convergence is \emph{exponentially fast} if there
exist $0<\sigma<1$ and $D>0$ such that for all $n\ge 1$ we have
$|x_n-x|\le D\sigma^n$.

\begin{defn}[Genericity]
  Let $S\subseteq \mathcal{W} \subseteq F$. We say that $S$ is
  \emph{exponentially $\mathcal{W}$-generic} if
\[
\lim_{n\to\infty} \frac{\gamma(n,S)}{\gamma(n,\mathcal{W})}=1
\]
and the convergence is exponentially fast. The complement in
$\mathcal{W}$ of an exponentially $\mathcal{W}$-generic set is called
\emph{exponentially $\mathcal{W}$-negligible}.
\end{defn}
In practice we are only interested in the cases $\mathcal{W} =F$ and
$\mathcal{W} =\mathcal{CR}$, the set of all cyclically reduced words
in $F$.  By definition $S\subseteq F$ is exponentially $F$-generic if
and only if $\lim_{n\to\infty} P_n(S)=1$ with exponentially fast
convergence in this limit. Similarly $S\subseteq \mathcal{CR}$ is
exponentially $\mathcal{CR}$-generic iff $\lim_{n\to\infty} P_n'(S)=1$
with exponentially fast convergence. Here there is a simple criterion
of being exponentially negligible~\cite{KSS} in $F$ and
$\mathcal{CR}$:
\begin{lem}\label{genericity}
  Let $\mathcal{W} =F$ or $\mathcal{W} =\mathcal{CR}$. Then for a
  subset $S\subseteq \mathcal{W} $ the following are equivalent:
\begin{enumerate}

\item The set $S$ is exponentially $\mathcal{W}$-negligible.

\item We have
\[ \frac{\gamma(n,S)}{(2k-1)^n} \to_{n\to\infty} 0
\text{ exponentially fast, }
\]

\item We have
\[ \frac{\rho(n,S)}{(2k-1)^n} \to_{n\to\infty} 0
\text{ exponentially fast, }
\]

\item We have
\[
\limsup_{n\to\infty} \sqrt[n]{\rho(n,S)} < 2k-1.
\]

\item We have
\[
\limsup_{n\to\infty} \sqrt[n]{\gamma(n,S)} < 2k-1.
\]
\end{enumerate}
\end{lem}

\begin{prop}\label{ld'}
  Let $\epsilon>0$ and let $m>0$ be an integer. Then the set

\begin{align*}
  W(m,\epsilon):=\{& w\in F : \text{ for any } u\ne 1 \text{ with }
  |u|=m \text{ we have }\\
  & \big|f_w(u) - \frac{1}{2k (2k-1)^{m-1}}\big|<\epsilon\}
\end{align*}
is exponentially $F$-generic.
\end{prop}

\begin{proof}
  This is a straightforward corollary of Large Deviation
  Theory~\cite{DZ} applied to the finite state Markov chain generating
  the freely reduced words in $F$. We refer the reader to~\cite{KSS}
  for a more detailed discussion about large Deviation Theory and how
  it works in this particular case.
\end{proof}

It is not hard to deduce the following from Proposition~\ref{ld'}.

\begin{prop}\label{ld}
  Let $\epsilon>0$ and let $m>0$ be an integer. Then the set

\begin{align*}
  C(m,\epsilon ):=\{& w\in \mathcal{CR} : \text{ for any } u\ne 1
  \text{ with } |u|=m \text{ and for the cyclic word } (w) \\ &\text{
    we have } \big|f_{(w)}(u) -
  \frac{1}{2k(2k-1)^{m-1}}\big|<\epsilon\}
\end{align*}
is exponentially $\mathcal{CR}$-generic.

\end{prop}

We now give the definition of an ``approximate'' stretching factor,
which will later be seen to be equivalent to the generic stretching
factor of an automorphism introduced earlier.

\begin{defn}
  Let $\phi:F\to Aut(X)$ be a free simplicial action without
  inversions of $F=F(a_1,\dots, a_k)$ on a simplicial tree $X$.

  We say that a number $\lambda\ge 0$ is a \emph{approximate
    stretching factor} of $\phi$ if for every $p\in X$ and for any
  $\epsilon>0$ the set
\[
\{w\in F: |\frac{|w|_p}{|w|}-\lambda|\le \epsilon\}
\]
is exponentially generic in $F$.

Similarly, we say that a number $\lambda\ge 0$ is a \emph{approximate
  conjugacy stretching factor} of $\phi$ if for any $\epsilon>0$ the
set
\[
\{w\in \mathcal{CR}: |\frac{||w||_X}{||w||}-\lambda|\le \epsilon\}
\]
is exponentially generic in $\mathcal{CR}$.

\end{defn}

\begin{prop}\label{equiv}
  Let $\phi:F\to Aut(X)$ be a free simplicial action of
  $F=F(a_1,\dots, a_k)$ on a simplicial tree $X$.

\begin{enumerate}
\item There is at most one approximate stretching factor for $\phi$.
\item There is at most one approximate conjugacy stretching factor for
  $\phi$.
\item If $\lambda$ is an approximate conjugacy stretching factor for
  $\phi$ then $\lambda$ is also an approximate stretching factor for
  $\phi$.
\item If $\lambda$ is an approximate stretching factor for $\phi$ then
  $\lambda$ is also an approximate conjugacy stretching factor for
  $\phi$.
\end{enumerate}
\end{prop}

\begin{proof}

  Parts (1) and (2) are obvious.

  We now establish (3).  Indeed, suppose that $\lambda$ is an
  approximate conjugacy stretching factor for $\phi$. Let $\epsilon>0$
  and let $S$ be the set of all cyclically reduced words $w$ such that
  $\big| \frac{||w||_X}{||w||}-\lambda\big|\ge \epsilon/2$.  Then $S$
  is exponentially $\mathcal{CR}$-negligible, so that

\[
\frac{\gamma(n,S)}{(2k-1)^n}\to_{n\to\infty} 0
\]
exponentially fast.  Let $p\in X$.  Put $M=\max\{ |a_i|_p: 1\le i\le
k\}$. Let $N>0$ be an integer such that for any cyclically reduced
word $u$ we have $\big||u|_p-||u||_X \big|\le N$.

Let $S'$ be the set of all freely reduced words $w$ in $F$ that differ
from an element of $S$ in at most the last letter. Then
$\gamma(n,S')\le 2k \gamma(n,S)$ and hence $S'$ is exponentially
$F$-negligible by Lemma~\ref{genericity}.

Suppose $w\in F-S'$ is such that
$\frac{N+2M}{|w|}<\frac{\epsilon}{2}$. Let $u$ be a cyclically reduced
word obtained from $w$ by changing at most the last letter. Then
$|u|=|w|$ and $u\in \mathcal{CR}-S$.

Thus $d_A(w,u)\le 2$ and hence $d_X(wp, up)\le 2M$. Thus
$\big||u|_p-|w|_p\big|\le 2M$. Also $\big|||u||_X-|u|_p\big|\le N$.
Therefore $\big|||u||_X-|w|_p\big|\le N+2M$.  Since $u\in
\mathcal{CR}-S$, we have $\big| ||u||_X -\lambda ||u|| \big|< \epsilon
||u||$. Since $||u||=|u|=|w|$, we have:

\begin{gather*}
  \big| |w|_p -\lambda |w| \big| < \epsilon |w|+N+2M\\
  \left| \frac{|w|_p}{|w|} -\lambda \right| < \frac{\epsilon}{2}
  +\frac{N+2M}{|w|}<\epsilon.
\end{gather*}

The set $\{w\in F : \frac{N+2M}{|w|}\ge \frac{\epsilon}{2}\}$ is
finite.  Hence $S'\cup \{w\in F : \frac{N+2M}{|w|}\ge
\frac{\epsilon}{2}\}$ is exponentially $F$-negligible and assertion
(3) holds.

The proof of (4) is similar to that of (3) and we leave the details to
the reader.
\end{proof}

\begin{thm}\label{thm:interm}
  Let $F=F(a_1,\dots, a_k)$ and let $\phi:F\to Aut(X)$ be a free
  simplicial action of $F=F(a_1,\dots, a_k)$ on a simplicial tree $X$.

  Then the generic stretching factor $\lambda(\phi)$ is also an
  approximate conjugacy stretching factor (and thus by
  Proposition~\ref{equiv} an approximate stretching factor).
\end{thm}

\begin{proof}
  The proof is very similar to that of Theorem~\ref{main}.  Since the
  minimal $F$-invariant subtree of $X$ contains the axes of all the
  nontrivial elements of $F$, we may again assume that the action of
  $F$ on $X$ is minimal.

  Choose a vertex $x\in X$.  Recall, that, using the notations from
  the proof of Theorem~\ref{main}, for any $w\in F$
\[
\big|\sum_{a\in Z_1} \sum_{|u|_A=M} c(u,a) f_{w}(u) + \sum_{z\in
Z_2} \sum_{|u|_A=M} r_z c(u,z) f_{w}(u)-\frac{|g|_p}{n}\big|\le
\frac{L'}{n} \tag{!}.
\]

It follows from Lemma~\ref{lem:rewr} and Lemma~\ref{lem:compute}
that if $w\in F$ is cyclically reduced over $B_T$ then $\big|
||w||_X-|w|_x\big|\le N$ where $N=N(x)>0$ is some constant
independent of $w$.  On the other hand by the Bounded Cancellation
Lemma ( Lemma~\ref{bnd}) there exists a constant $N'>0$ such that
for any cyclically reduced word $w$ over $A$, we have $\big|
||w||_{B_T}-|w|_{B_T}\big|\le N'$. Hence for any cyclically
reduced word $w$ over $A$ we have $\big| ||w||_X-|w|_x\big|\le
N''$ where $N''=N''(x)>0$ is some constant independent of $w$.

Therefore for any cyclically reduced $w\in F$ over $A$ with
$||w||=n$
\[
\left| \sum_{a\in Z_1} \sum_{|u|_A=M} c(u,a) f_{w}(u) + \sum_{z\in
    Z_2} \sum_{|u|_A=M} r_z c(u,z) f_{w}(u)-\frac{||w||_X}{n} \right|
\le \frac{L'+N}{n} \tag{*}.
\]

Let $\epsilon>0$ We know that the set

\begin{align*}
  C(M,\epsilon):=\{& w\in \mathcal{CR} : \text{ for any } u\ne 1
  \text{ with } |u|=M \text{ and for the cyclic word } (w) \\ &\text{
    we have } \left|f_{(w)}(u) - \frac{1}{2k(2k-1)^{M-1}} \right|
  <\epsilon\}
\end{align*}
is exponentially $\mathcal{CR}$-generic.

Hence (*) implies that for any $w\in C(M,\epsilon)$

\[
\left| \frac{1}{2k(2k-1)^{M-1}} \left( \sum_{a\in Z_1} \sum_{|u|_A=M}
    c(u,a) + \sum_{z\in Z_2} \sum_{|u|_A=M} r_z c(u,z) \right)
  -\frac{||w||_X}{n} \right| \le \frac{N_1}{n}+N_1\epsilon.
\]
for some constant $N_1>0$ independent of $w$ and $\epsilon$.

Thus by definition the number

\[
\frac{1}{2k(2k-1)^{M-1}} \left( \sum_{a\in Z_1} \sum_{|u|_A=M} c(u,a)
  + \sum_{z\in Z_2} \sum_{|u|_A=M} r_z c(u,z) \right)
\]
is an approximate conjugacy stretching factor for $\phi$. In the proof
of Theorem~\ref{main} we obtained the same formula for
$\lambda(\phi)$.
\end{proof}

\begin{lem}\label{aux}

  Let $F=F(a_1,\dots, a_k)$ and let $\phi:F\to Aut(X)$ be a free
  simplicial action of $F=F(a_1,\dots, a_k)$ on a simplicial tree $X$.
  Let $\mu\ge 0$.

  Suppose there exists an exponentially $\mathcal{CR}$-generic set $S$
  such that for any $w\in S$ \[ \frac{||w||_X}{||w||}\ge \mu.
\]

Then $\lambda(\phi)\ge \mu$.
\end{lem}

\begin{proof}

  Suppose, on the contrary, that $\lambda(\phi)< \mu$. Choose
  $\epsilon>0$ such that $\lambda(\phi)+\epsilon<\mu$.

  Then there is an exponentially $\mathcal{CR}$-generic set $R$ of
  cyclically reduced words such that for any $w\in R$

\[
\frac{||w||_X}{||w||} \le \lambda+\epsilon.
\]
The intersection $S\cap R$ is exponentially $\mathcal{CR}$-generic and
hence nonempty. Take $w\in S\cap R$.

Then

\[
\mu \le \frac{||w||_X}{||w||}\le \lambda+\epsilon<\mu,
\]
yielding a contradiction.
\end{proof}

\section{Whitehead's Peak Reduction and rigidity of free group automorphisms}\label{sect:wh}

We need to recall some definitions related to Whitehead's algorithm
for solving the automorphic equivalence problem in a free group. We
refer the reader to~\cite{LS,Wh} for a detailed exposition.

\begin{defn}[Whitehead automorphisms]\label{defn:moves}
  A \emph{Whitehead automorphism} of $F$ is an automorphism $\tau$ of
  $F$ of one of the following two types:

  (1) There is a permutation $t$ of $\widehat{A}$ such that
  $\tau|_{\widehat{A}}=t$. In this case $\tau$ is called a
  \emph{relabeling automorphism} or a \emph{Whitehead automorphism of
    the first kind}.

  (2) There is an element $a\in \widehat{A}$, the \emph{multiplier},
  such that for any $x\in \widehat{A}$
\[
\tau(x)\in \{x, xa, a^{-1}x, a^{-1}xa\}.
\]

In this case we say that $\tau$ is a \emph{Whitehead automorphism of
  the second kind}. (Note that we always have $\tau(a)=a$ in this case
since $\tau$ is an automorphism of $F$.) To every such $\tau$ we
associate a pair $(S,a)$ where $a$ is as above and $S$ consists of all
those elements of $\widehat{A}$, including $a$ but excluding $a^{-1}$,
such that $\tau(x)\in\{xa, a^{-1}xa\}$.  We will say that $(S,a)$ is
the \emph{characteristic pair} of $\tau$.
\end{defn}
Note that for any $a\in \widehat{A}$ the inner automorphism $ad(a)$ is
a Whitehead automorphism of the second kind.

The following important result of Whitehead is known as the ``peak
reduction lemma'':

\begin{prop}\label{peak_red}
  Let $u,v$ be cyclic words with $||u||\le ||v||$ and let $\phi\in
  Aut(F)$ be such that $\phi(u)=v$. Then we can write $\phi$ as a
  product of Whitehead moves

\[
\phi=\tau_p\dots \tau_1
\]
so that for each $i=1,\dots, p$
\[
||\tau_i\dots \tau_1(u)|| \le ||v||.
\]
Moreover, if $||u||<||v||$ then the above inequalities are strict for
all $i<p$.
\end{prop}

\begin{defn}[Weighted Whitehead graph]

  Let $w$ be a nontrivial cyclically reduced word in
  ${\widehat{A}}^*$.  The \emph{weighted Whitehead graph $\Gamma_w$ of
    $w$} is defined as follows.  Let $(w)$ be the cyclic word defined
  by $w$. The vertex set of $\Gamma_w$ is $\widehat{A}$. For every
  $x,y\in \widehat{A}$ such that $x\ne y^{-1}$ there is an undirected
  edge in $\Gamma_w$ from $x^{-1}$ to $y$ labeled by the sum $\hat
  n_{w}(xy):=n_{(w)}(xy)+ n_{(w)}(y^{-1}x^{-1})$.
\end{defn}

There are $k(2k-1)$ undirected edges in $\Gamma_w$. Edges may have
label zero, but there are no edges from $a$ to $a$ for $a\in
\widehat{A}$.  It is easy to see that we have $\Gamma_{w}=\Gamma_v$
for any cyclic permutation $v$ of $w$ or $w^{-1}$.

\begin{conv}
  Let $w$ be a fixed nontrivial cyclically reduced word.  For two
  subsets $X,Y\subseteq \widehat{A}$ we denote by $X.Y$ the sum of all
  edge-labels in the weighted Whitehead graph $\Gamma_w$ of $w$ of
  edges from elements of $X$ to elements of $Y$. Thus for $x\in
  \widehat{A}$ the number $x.\widehat{A}$ is equal to
  $n_w(x)+n_w(x^{-1})$, the total number of occurrences of $x^{\pm 1}$
  in $w$.
\end{conv}

The next lemma, which is Proposition~4.16 of Ch.~I in \cite{LS}, gives
an explicit formula for the difference of the lengths of $w$ and
$\tau(w)$, where $\tau$ is a Whitehead automorphism.

\begin{lem}\label{lem:LS}
  Let $w$ be a nontrivial cyclically reduced word and let $\tau$ be a
  Whitehead automorphism of the second kind with the characteristic
  pair $(S,a)$. Let $S'=\widehat{A}-S$. Then
 \[  ||\tau(w)||-||w||=S.S'-a.\widehat{A}. \]
\end{lem}

The following important notion was introduced by Kapovich, Schupp and
Shpilrain in \cite{KSS}.

\begin{defn}[Strict Minimality]
  A nontrivial cyclically reduced word $w$ in $F$ is \emph{strictly
    minimal} if for every non-inner Whitehead automorphism $\tau$ of
  $F$ of the second kind we have
\[
||\tau(w)||>||w||.
\]
The set of all strictly minimal elements in $F$ is denoted $SM$.
\end{defn}

An immediate consequence of the Peak Reduction Lemma is:

\begin{prop}\cite{KSS}
  Let $w \in F$ be a cyclically reduced strictly minimal element. Then
  $w$ is of minimal length in its $Aut(F)$-orbit and for any $\phi\in
  Aut(F)$ we have:

\[
|w|=||w||\le ||\phi(w)||\le |\phi(w)|.
\]
\end{prop}

\begin{thm}\label{B}
  Put $c_0:=1+\frac{2k-3}{4k^2-2k}$. There exists an exponentially
  $F$-generic set $W\subseteq F$ with the following property.

  For any $\phi\in Aut(F)$ the following conditions are equivalent:

\begin{enumerate}

\item The automorphism $\phi$ is simple.

\item We have $\lambda(\phi)=1.$

\item We have $\lambda(\phi)< 1+\frac{2k-3}{2k^2-k}.$

\item For some $w\in W$ we have $||\phi(w)||=||w||$.

\item For every $w\in W$ we have $||\phi(w)||=||w||$.

\item For some $w\in W$ we have $||\phi(w)||\le c_0 ||w||$.

\item For every $w\in W$ we have $||\phi(w)||\le c_0 ||w||$.

\end{enumerate}

\end{thm}

\begin{proof}

  It is obvious that (1) implies (2) and that (2) implies (3).

  We will now prove that (3) implies (1).

  Let $\phi\in Aut(F)$.

  Let $\epsilon>0$ be arbitrary. Put $T(\epsilon)$ be the set of all
  cyclically reduced words $w$ such that:

\begin{itemize}
\item For any $x\in \widehat{A} $ $|f_w(x)-\frac{1}{2k}|\le
  \epsilon/2$.

\item For any $x,y\in \widehat{A} $ with $x\ne y^{-1}$
  $|f_w(xy)-\frac{1}{2k(2k-1)}|\le \epsilon/2$
\end{itemize}

Then $T(\epsilon)$ is exponentially $\mathcal{CR}$-generic~\cite{KSS}.
Moreover, every $w\in T(\epsilon)$ is strictly minimal~\cite{KSS},
provided that $\epsilon<\frac{2k-3}{k(2k - 1) (4k - 3)}$.

Suppose now that $\epsilon<\epsilon_0:=\frac{2k-3}{k(2k - 1) (4k -
  3)}$.  Choose an arbitrary element $w\in T(\epsilon)$ and denote
$n=||w||$.

By strict minimality of $w$ we have $||w||\le ||\phi(w)||$.  Moreover,
by Proposition~\ref{peak_red} (Peak Reduction Lemma) there is a
decomposition $\phi=\tau_p \tau_{p-1}\dots \tau_1$ such that each
$\tau_i$ is a Whitehead move and for each $i=1,\dots, p-1$
\[
||\tau_i \tau_{i-1}\dots \tau_1(w)||\le ||\phi(w)||
\]
with strict inequalities unless $||w||=||\phi(w)||$.

Suppose first that $||w||=||\phi(w)||$. Then all inequalities above
are equalities and by strict minimality of $w$ each $\tau_i$ is either
inner or a relabeling automorphism. This implies that $\phi=\alpha
\tau$ where $\alpha$ is inner and $\tau$ is a relabeling automorphism
and that $\lambda(\phi)=1$.

Suppose now that $||w||< ||\phi(w)||$. Then the preceding argument
shows that in fact for any $z\in T(\epsilon)$ we have $||z||<
||\phi(z)||$ (since otherwise $\phi$ would be simple and so
$||w||=||\phi(w)||$).

Denote $z_0=z$ and $z_i=\tau_i \tau_{i-1}\dots \tau_1(z)$ for $0<i\le
p$. Thus $z_p=\phi(z)$. Since $||z||< ||\phi(z)||$, there is some $i,
1\le i\le p$ such that $\tau_i$ is a non-inner Whitehead move of the
second kind. Let $j$ be the smallest $i$ with this property. Then all
$\tau_{i}$ with $i<j$ are either inner or relabeling automorphisms,
$||z||=||z_{i}||$ and $z_{i}\in T(\epsilon)$. In particular,
$z_{j-1}\in T(\epsilon)$ and $z_{j-1}$ is strictly minimal.

Thus \[n=||z||=||z_{j-1}||\le
||z_j||=||\tau_j(z_{j-1})||<||\phi(z)||.\]

Let $(S,a)$ be the characteristic pair of $\tau_j$ and let
$S'=\widehat{A}-S$.  Since $\tau_j$ is non-inner, we have both $|S|
\ge 2$, and $|S'| \ge 2$. Hence $|S|\ |S'| \ge 2(2k-2)$ and there are
at least $2(2k-2)$ edges between $S$ and $S'$ in the weighted
Whitehead graph of $z_{j-1}$. Recall that $a. \widehat{A}$ is the
total number of occurrences of $a^{\pm 1}$ in $z$.

By Lemma~\ref{lem:LS}, we have
$||\tau_j(z_{j-1})||-||z||=S.S'-a.\widehat{A}$.

By assumption on $z_{j-1}$ we have $a.\widehat{A}\le
n(\frac{1}{k}+\epsilon)$ and so

\[
||\tau_j(z_{j-1})||-||z_{j-1}||=S.S'-a.\widehat{A} \ge 2n(2k-2) \left(
  \frac{1}{k(2k-1)}-\epsilon \right)- n\left( \frac{1}{k}+\epsilon
\right).
\]

Hence

\[
||\phi(z)||\ge ||z_j||\ge n+ 2n(2k-2) \left(
  \frac{1}{k(2k-1)}-\epsilon \right) - n \left( \frac{1}{k}+\epsilon
\right)
\]

and so, since $n=||z||$,

\[
\frac{||\phi(z)||}{||z||}\ge 1+ (4k-4) \left(
  \frac{1}{k(2k-1)}-\epsilon \right)- \left( \frac{1}{k}+\epsilon
\right)
\]

Note that the above inequality holds for any element $z\in
T(\epsilon)$.

Since $T(\epsilon)$ is exponentially $\mathcal{CR}$-generic, this
implies by Lemma~\ref{aux} that
\[
\lambda(\phi)\ge 1+
(4k-4)(\frac{1}{k(2k-1)}-\epsilon)-(\frac{1}{k}+\epsilon).
\]
Since $0<\epsilon<\epsilon_0$ was arbitrary, it follows that

\[
\lambda(\phi)\ge 1+
(4k-4)\frac{1}{k(2k-1)}-\frac{1}{k}=1+\frac{2k-3}{2k^2-k}>1.
\]

This proves that (3) implies (1), so that (1), (2) and (3) are
equivalent.

Choose $0<\epsilon_1<\epsilon_0$ such that

\[
1+ (4k-4) \left( \frac{1}{k(2k-1)}-\epsilon_1 \right) - \left(
  \frac{1}{k}+\epsilon_1 \right) <c_0=1+\frac{2k-3}{4k^2-2k} .
\]

Put $W=T(\epsilon_1)$. The above argument shows that if for some $w\in
W$ we have
\[
\frac{||\phi(z)||}{||z||}< 1+ (4k-4) \left(
  \frac{1}{k(2k-1)}-\epsilon_1 \right) - \left( \frac{1}{k}+\epsilon_1
\right)
\]
then $\phi$ is simple.

With this $W$ we have proved that (5) implies (1). It is obvious that
(1) implies (4)-(7) and that each of (4), (5), (7) implies (6).  Thus
statements (1), (4), (5), (6), (7) are equivalent.

We already know that (1), (2) and (3) are equivalent. This completes
the proof of the theorem.
\end{proof}

The following statement, together with Theorem~\ref{B}, implies
Theorem~\ref{rigid} from the Introduction.

\begin{cor}\label{rough}
  Let $F=F(a_1,\dots,a_k)$, where $k\ge 2$, and $d$ be the word metric
  on $F$ corresponding to the generating set $A=\{a_1,\dots, a_k\}$.
  Let $\phi\in Aut(F)$.  Then the following conditions are equivalent

\begin{enumerate}
\item The automorphism $\phi$ is simple.
\item The map $\phi:(F,d)\to (F,d)$ is a rough isometry.
\item The map $\phi:(F,d)\to (F,d)$ is a rough similarity.
\end{enumerate}
\end{cor}

\begin{proof}
  It is obvious that (1) implies (2) and that (2) implies (3).

  We will now show that (3) implies (1). Suppose that $\phi$ is a
  rough similarity, so that there exist $\lambda >0$ and $D>0$ such
  that for any $w\in F$

\[
\lambda |w|-D \le |\phi(w)|\le \lambda |w|+D.
\]

Then obviously $\lambda=\lambda(\phi)$.  By Theorem~\ref{B} either
$\phi$ is simple or $\lambda(\phi)>1$.

Assume the latter. Put $\lambda_0=\frac{1+\lambda}{2}$. Thus
$1<\lambda_0<\lambda$.

Consider the ball $B_n$ of radius $n$ in $F$, where $n>>1$.  For any
$w\in F$ with $|w|\ge n/\lambda_0$ we have

\[
|\phi(w)|\ge \lambda |w|-D \ge \lambda n/\lambda_0 -D> n,
\]
so that $\phi(w)\not\in B_n$.

Thus only the elements of length $\le n/\lambda_0$ may be potentially
taken to $B_n$ by $\phi$. The number of such elements is smaller than
$\# B_n$ since $\lambda_0>1$ and $n/\lambda_0 <n$.  This contradicts
the fact that $\phi:(F,d)\to (F,d)$ is a bijection.  Therefore $\phi$
is simple, as required.
\end{proof}

The following is Theorem~\ref{rigid1} from the Introduction:

\begin{thm}\label{thm:jump-simpl}
  Let $F=F(a_1,\dots, a_k)$ where $k\ge 2$. Let $\phi: F\to X$ be a
  free minimal action on $F$ on a simplicial tree $X$ without
  inversions.

  Then exactly one of the following occurs:

\begin{enumerate}
\item There is a simple automorphism $\alpha$ of $F$ such that $X$ is
  $\phi\circ \alpha$-equivariantly isomorphic to the Cayley graph of
  $F$ with respect to $\{a_1,\dots, a_k\}$. In this case
  $\lambda(\phi)=1$.
\item We have $\lambda(\phi)\ge 1+\frac{1}{k(2k-1)}$.
\end{enumerate}
\end{thm}

\begin{proof}
  Let $\Gamma=X/F$ and let $T$ be a maximal tree in $\Gamma$ and let
  $B=\{b_1,\dots, b_k\}$ be the geometric basis of $F$ corresponding
  to $T$.  Let $\psi\in Aut(F)$ be defined by $\alpha(b_i)=a_i$ for
  $i=1,\dots, k$.

  Note that because of Lemma~\ref{lem:compute} for any cyclic word $w$
  over $B$ we have $||w||_X\ge ||w||_B$.  Suppose first that $\alpha$
  is not a simple automorphism. Then $\lambda(\alpha)\ge
  1+\frac{2k-3}{k(2k-1)}$.

  Hence for every $\epsilon>0$ there exists an exponentially
  $\mathcal{CR}$-generic set $R(\epsilon)\subseteq \mathcal{CR}$ such
  that for any $w\in R(\epsilon)$

\[
\frac{||w||_B}{||w||_A}=\frac{||\alpha(w)||_A}{||w||_A}\ge
1+\frac{2k-3}{k(2k-1)}-\epsilon.
\]
Since $||w||_X\ge ||w||_B$, it follows that
\[
\frac{||w||_X}{||w||_A}\ge 1+\frac{2k-3}{k(2k-1)}-\epsilon.
\]
Since $\epsilon>0$ was arbitrary, it follows by Lemma~\ref{aux} that
\[\lambda(\phi)\ge 1+\frac{2k-3}{k(2k-1)}\ge 1+\frac{1}{k(2k-1)},\] as
required.

Suppose now that $\alpha$ is a simple automorphism.  We will assume
that $\alpha=Id$, and it will be easily seen that the general case is
similar.

If $\Gamma$ is a wedge of $k$ loop-edges then the statement of the
theorem holds. Suppose $\Gamma$ is not of this form. Then there exist
edges $e,e'\in E(\Gamma-T)$, $e'\ne e^{-1}$, such that $[t(e),
o(e')]_T$ has length at least $1$. Let $z$ be the freely reduced word
of length $2$ in $B$ corresponding to the sequence $ee'$. Let
$\epsilon>0$ be arbitrary. Let $C(2,\epsilon/2)\subseteq \mathcal{CR}$
be defined as in Proposition~\ref{ld}. Thus $C(2,\epsilon/2)$ consists
of all cyclically reduced words $w'$ such that for the cyclic word
$w=(w')$ and for every freely reduced word $xy$ in $A$
\[
\big| f_w(xy)-\frac{1}{2k(2k-1)}\big|\le \epsilon/2.
\]
Then $C(2,\epsilon/2)$ is exponentially $\mathcal{CR}$-generic.  Let
$w'\in C(2,\epsilon/2)$ be arbitrary and let $w=(w')$. Note that
$||w'||_A=||w'||_B=||w||_A=||w||_B$ and $||w||_X=||w'||_X$.

Then

\[
||w||_X\ge ||w||_B+n_w(z)+n_w(z^{-1})=||w||_A+n_w(z)+n_w(z^{-1})
\]
and so
\[
\frac{||w'||_X}{||w'||_A}=\frac{||w||_X}{||w||_A}\ge
1+f_w(z)+f_w(z^{-1})\ge 1+\frac{1}{k(2k-1)}-\epsilon.
\]
Since $\epsilon>0$ was arbitrary, Lemma~\ref{aux} implies that
$\lambda(\phi)\ge 1+\frac{1}{k(2k-1)}$, as required.

\end{proof}

\section{Application to the geometry of automorphisms}\label{sect:appl}

\begin{defn}
  Let $F=F(a_1,\dots, a_k)$. An automorphism $\phi$ of $F$ is said to
  be \emph{$(s,m)$-hyperbolic}, where $s>1$ and $m\ge 1$ is an
  integer, if for every nontrivial cyclic word $w$ we have
\[
s ||w|| \le \max \{ ||\phi^m(w)||, ||\phi^{-m}(w)||\}.
\]
An automorphism is \emph{hyperbolic} if it is $(s,m)$-hyperbolic for
some $s>1, m\ge 1$.
\end{defn}

The following lemma is an easy consequence of the above definition:

\begin{lem}\label{gr}
  Let $\phi\in Aut(F)$ be $(s,m)$-hyperbolic and let $w$ be a cyclic
  word of minimal length in its $\langle \phi \rangle$-orbit.  Then
  for any $n\ge 2$ we have
\[
||\phi^{mn}(w)||\ge s^{n-1}||w||.
\]
\end{lem}

\begin{proof}

  By definition of hyperbolicity of $\phi$ we have
\[
||\phi^{-m}(u)||\le ||u|| \Rightarrow s ||u||\le ||\phi^m(u)||
\tag{!}.
\]

Note that by the choice of $w$ we have $||w||\le ||\phi^{m}(w)||$.
Hence applying (!) with $u=\phi^m(w)$ we get $s ||\phi^m(w)||\le
||\phi^{2m}(w)||$. Then, using (!), by induction on $n$ we see
that for any $n\ge 1$
\[
||\phi^{m(n+1)}(w)||=||\phi^{mn+m}(w)||\ge s ||\phi^m(w)||.
\]
This in turn implies that for any $n\ge 1$
\[
||\phi^{m(n+1)}(w)||=||\phi^{mn+m}(w)||\ge s^n ||\phi^m(w)||\ge
s^{n-1} ||w||.
\]
This proves Lemma~\ref{gr}.
\end{proof}

The following is Theorem~\ref{hyp} from the Introduction:

\begin{thm}\label{hyperb}
  Let $\phi$ be an $(s,m)$-hyperbolic automorphism of $F$. Then
\[
\liminf_{n\to\infty}\sqrt[n]{\lambda(\phi^n)}\ge \sqrt[m]{s}>1.
\]
\end{thm}

\begin{proof}
  Let $t\ge 2$ be an arbitrary integer. Let $w\in SM$ be a strictly
  minimal element. Since $w$ is minimal in its $Aut(F)$-orbit, it is
  also minimal in its $\langle \phi \rangle$-orbit. Therefore by
  Lemma~\ref{gr}
\[
||\phi^{tm}(w)||\ge s^{t-1} ||w|| \text{ and }
\frac{||\phi^{tm}(w)||}{||w||}\ge s^{t-1}.
\]
Since $SM$ is exponentially $\mathcal{CR}$-generic, Lemma~\ref{aux}
implies that $\lambda(\phi^{tm})\ge s^{t-1}$.

Moreover, there is $D>0$ such that for any cyclically reduced word $u$
we have

\[
||\phi^i(u)||\ge D ||u||, \text{ for all } 0\le i <m.
\]

Let $n\ge 2m$ be an integer. Then we can write $n$ as $n=mt+i$
where $t\ge 2$ and $0\le i<m$. For any $w\in SM$ we have

\[
||\phi^{n}(w)||=||\phi^{mt+i}(w)||\ge D ||\phi^{mt}(w)||\ge D s^{t-1}
||w||
\]
and hence
\[
\frac{||\phi^{tm}(w)||}{||w||}\ge Ds^{t-1}.
\]
Since $SM$ is exponentially $\mathcal{CR}$-generic, Lemma~\ref{aux}
again implies that for any $n\ge 2m$
\[
\lambda(\phi^n)\ge Ds^{t-1}= \frac{D}{s} s^{\lfloor n/m\rfloor}.
\]
This implies
\[
\liminf_{n\to\infty}\sqrt[n]{\lambda(\phi^n)}\ge s^{1/m}>1,
\]
as claimed.
\end{proof}

\setcounter{theore}{6}

We can now prove Corollary~\ref{cor:flux} from the Introduction:

\begin{cor}
  Let $F=F(a_1,\dots, a_k)$ be a free group of rank $k\ge 2$, equipped
  with the standard metric.

  Then for any $\phi\in Aut(F)$ we have:

\[
flux(\phi)=\begin{cases}
  0, \text{ if } \phi \text{ is a    relabeling automorphism}\\
  1, \text{ otherwise.}
\end{cases}
\]
\end{cor}
\begin{proof}
  Let $\lambda=\lambda(\phi)$ be the generic stretching factor.

  Suppose first that $\lambda>1$. Then the set
\[
T:=\{ w \in F : \frac{|\phi(w)|}{|w|}> \frac{\lambda+1}{2}
\]
is exponentially $F$-generic.

Let $B(n)$ be the ball of radius $n$ in $F$ and let $w \in B(n)\cap T$
be such that $2n/(\lambda+1)\le |w|\le n$. Then

\[
|\phi(w)|> |w|\frac{\lambda+1}{2}\ge \frac{2n}{\lambda+1}
\frac{\lambda+1}{2}=n
\]
Hence for each $w\in [B(n)\cap T]-B(2n/(\lambda+1))$ we have
$|\phi(w)|>n$. The size of $B(2n/(\lambda+1))$ is exponentially
smaller than that of $B(n)$ since $2/(\lambda+1)<1$. Hence by
exponential genericity of $T$
\[
\frac{\#[B(n)\cap T]-\#B(2n/(\lambda+1))}{\#
  B(n)}\longrightarrow_{n\to \infty} 1 \text{ exponentially fast}.
\]
Hence \[\lim_{n\to\infty}\frac{flux_{\phi}(n)}{\# B(n)}=1\] and
therefore $flux(\phi)=1$.

Suppose now that $\lambda(\phi)=1$. By Theorem~\ref{B} this implies
that $\phi=\alpha \tau$ where $\alpha$ is inner and $\tau$ is a
relabeling automorphism.

If $\alpha=1$, then obviously $flux(\phi)=0$. Suppose now that
$\alpha$ is nontrivial. Since $\tau$ acts as a permutation on each
ball and each sphere in $F$, we can assume that $\tau=1$ and
$\phi=\alpha$. Thus there is $u\in F, u\ne 1$ such that for every
$w\in F$ $\phi(w)=uwu^{-1}$. There are $\ge c_1 (2k-1)^n$ elements $f$
with $|w|=n$ such that the product $uwu^{-1}$ is freely reduced as
written, where $c_1>0$ is a constant independent of $n$ and $u$. For
each such element we have $|\phi(w)|>|w|$. Hence there is a constant
$c_2\in (0,1)$ independent of $n$ and $u$ such that for any $n>0$

\[
1\ge \frac{flux_{\phi}(n)}{\# B(n)}\ge c_2>0.
\]
Hence
\[
1\ge flux(\phi)= \lim_{n\to\infty} \sqrt[n]{\frac{flux_{\phi}(n)}{\#
    B(n)}} \ge \lim_{n\to\infty}\sqrt[n]{c_2}=1.
\]
Thus $flux(\phi)=1$ and the proof is complete.

\end{proof}

\section{Random elements in regular languages}\label{sect:lang}

The most reasonable way of choosing a ``random'' element in the
regular language $L$ is via a random walk in the transition graph of
an automaton $M$ accepting $L$.  It turns out that the natural model
of computation here is that of a \emph{non-deterministic finite
  automaton} or NDFA. Such an automaton $M$ over an alphabet $A$ with
\emph{state set} $Q$ is specified by a finite directed graph
$\Gamma(M)$.  The vertex set of $\Gamma(M)$ is the set $Q$ of states
of $M$ and $Q$ comes equipped with a distinguished nonempty subset $I$
of \emph{initial} or \emph{start} states.  The directed edges of
$\Gamma(M)$ are labelled by elements of $A$ and these edges are
treated as transitions of $M$. If $q\in Q$ is a state and $a\in A$ is
a letter, we allow multiple edges labelled $a$ with origin $q$ and we
also allow the case when there are no such edges. Nondeterministic
automata are thus by their nature "partial". There is a distinguished
subset of $Q$ of \emph{accepting} states. A word $w$ over $A$ is said
to be \emph{accepted} by $M$ if there exists a directed path with
label $w$ in $\Gamma(M)$ from some initial state to an accepting
state.  The \emph{language}, $L(M)$, accepted by $M$ is the collection
of all words accepted by $M$.

We will also use directed graph $\Gamma_1(M)$ defined as follows.  The
vertex set of of $\Gamma_1(M)$ is the set of directed edges
$E(\Gamma(M))$ of $M$.  If $e_1,e_2\in E(\Gamma(M))$ the pair
$(e_1,e_2)$ defines a directed edge from $e_1$ to $e_2$ in
$\Gamma_1(M)$ if the terminus of $e_1$ is the origin of $e_2$, that
is, $e_1,e_2$ is a directed edge-path in $\Gamma(M)$.

\begin{defn}[Normal  Automaton]\label{defn:normal}
  Let $A$ be a finite alphabet.  A \emph{normal} automaton over a
  finite alphabet $A$ is a nondeterministic finite state automaton $M$
  over $A$ such that the following conditions hold:
\begin{itemize}
\item the automaton $M$ has a nonempty set of accept states;

\item the directed graph $\Gamma(M)$ has at least one edge;

\item the directed graph $\Gamma(M)$ is \emph{strongly connected},
  that is for any two states $q,q'$ of $M$ there exists a directed
  edge-path from $q$ to $q'$ in $\Gamma(M)$.
\end{itemize}
\end{defn}

The third condition in the above definition is the most important one
as it is responsible for the irreducibility of a Markov chain
naturally associated to a normal automaton:

\begin{defn}[Associated Markov chain]
  Let $M$ be a normal automaton over a finite alphabet $A$. We define
  an \emph{associated finite state Markov chain} $M'$ as follows. The
  set of states of $M'$ is the set $E$ of directed edges of
  $\Gamma(M)$.  If the origin of $f$ is not the terminus of $e$ we put
  the transition probability $p_{e,f}=0$. If the origin of $f$ is
  equal to the terminus of $e$ we put $p_{e,f}=1/m$, where $m$ is the
  total number of outgoing directed edges from the terminus of $e$.
\end{defn}

\begin{conv}
  Note that the sample space $\Omega$ for the Markov chain $M'$
  defined above consists of all semi-infinite directed edge-paths
\[
\omega=e_1,e_2,\dots, e_n, \dots
\]
in the graph $\Gamma(M)$.  Every such path has a label
\[
w(\omega)=a_1a_2\dots,
\]
that is a semi-infinite word over the alphabet $A$. We will denote
$w_n=w_n(\omega):=a_1\dots a_n$, the initial segment of length $n$ of
$w$. The set $\Omega$ comes equipped with the natural topology, where
we think of $\Omega$ as the union of boundaries of rooted trees
$(T_e)_e\in E$. The vertices of $T_e$ are finite edge-path in
$\Gamma(M)$ beginning with $e$. The Borel $\sigma$-algebra on $\Omega$
is generated by the following open-closed \emph{cylinder sets}
$Cyl(\gamma)$, where $\gamma$ is a nonempty finite edge-path in
$\Gamma(M)$:
\[
Cyl(\gamma):=\{\omega\in \Omega: p \text{ is the initial segment of }
\omega\}.
\]

If we put an initial probability distribution $\mu$ on $E$, this
defines a Borel probability measure $P_{\mu}$ on $\Omega$. This
measure is defined on the cylinder sets by the standard convolution
formula. If $\gamma=e_1,\dots, e_n$, where $n>1$, then

\[
P_{\mu}(Cyl(\gamma)):=\mu(e_1) p_{e_1,e_2} p_{e_2,e_3}\dots
p_{e_{n-1},e_n}.
\]
If $n=1$ then $P_{\mu}(Cyl(e)):=\mu(e)$.
\end{conv}

\begin{lem}
  Let $M$ be a normal automaton. Then the associated finite state
  Markov chain $M'$ is irreducible. In particular, there is a unique
  stationary initial probability distribution $\mu_0$ on the set of
  states $E$ of $M'$. This distribution has the property $\mu_0(e)>0$
  for each $e\in E$.
\end{lem}

\begin{proof}
  To show that $M'$ is irreducible we have two prove that for any two
  edges $e,f\in E$ there is $n>0$ such that the $n$-step transition
  probability $p^{(n)}_{e,f}>0$. Since $\Gamma(M)$ is strongly
  connected, there exists a directed edge-path $\gamma$ in $\Gamma(M)$
  from the terminus of $e$ to the origin of $f$. Then $e\gamma f$ is a
  directed edge-path in $\Gamma(M)$ that starts with $e$ and ends with
  $f$. Hence $\Gamma_1(M)$ is strongly connected and therefore $M'$ is
  irreducible.

  The irreducibility of $M'$ implies the existence and uniqueness of a
  positive stationary distribution $\mu_0$ on $E$, as required.
\end{proof}

If we fix an initial probability distribution $\mu$ on $E$, this
defines a probability measure $P_{\mu}$ on $\Omega$.

\begin{lem}\label{zero}
  Let $M$ be a normal automaton. Let $M'$ be the associated finite
  state Markov chain and let $\mu_0$ be the stationary initial
  distribution for $M'$. Let $Z \subseteq \Omega$ be a set such that
  $P_{\mu_0}(Z)=0$. Then for any other initial distribution $\mu$ on
  $E$ we have $P_{\mu}(Z)=0$.
\end{lem}

\begin{proof}
  Let $\mu$ and $\mu_0$ be as above. Put
\[
c:=\max\{ \frac{\mu(e)}{\mu_0(e)}: e\in E\}.
\]
Note that $0<c<\infty$ since $\mu_0(e)>0$ for each $e\in E$.  Consider
an arbitrary cylinder set $Cyl(\gamma)\subset \Omega$, where
$\gamma=e_1,e_2,\dots e_n$. From the definitions of $P_{\mu}$ and
$P_{\mu_0}$ we see that
\[
P_{\mu}(Cyl(\gamma))=\frac{\mu(e_1)}{\mu_0(e_1)}
P_{\mu_0}(Cyl(\gamma))\le c P_{\mu_0}(Cyl(\gamma)).
\]
Hence for an arbitrary Borel set $Z\subseteq \Omega$ we have
$P_{\mu}(Z)\le c P_{\mu_0}(Z)$. In particular, if $P_{\mu_0}(Z)=0$
then $P_{\mu}(Z)=0$.
\end{proof}

The previous two lemmas depend only on the automaton $M$ being normal.
Suppose now that $L = L(M)$.  For each state $q$ choose a shortest
path from $q$ to an accept state and let $u_q$ be the word in $A^\ast$
labelling that path. This is possible since $\Gamma(M)$ is strongly
connected and the set of accept states is nonempty by the assumption
on $M$.  Note that $u_q$ is the empty word if and only if $q$ is an
accept state. The lengths of $u_q$ are bounded above by some constant
depending on $M$.  For a finite walk $w_n$ denote $w_n'=w_nu_q$ where
$q$ is the state in which $w_n$ ends. Note that if $w_n$ begins in a
state from $I$ then $w_n'\in L$. Thus if $\mu$ is an distribution
supported on the set of edges in $E(M)$ with initial vertices from $I$
and $w_n$ is obtained by performing $n$ steps of the chain $M'$ with
initial distribution $\mu$, then $w_n'\in L$ can be thought of as a
"random" element of $L$.

We can now prove (a slight generalization of) Theorem~\ref{lang} from
the Introduction:

\begin{thm}\label{thm:H}
  Let $M$ be a normal automaton over the alphabet $A$ and let $L=L(M)$
  be the language accepted by $M$.

  Let $\phi: A^*\to G$ be a monoid homomorphism, where $G$ is a group
  with a left-invariant semi-metric $d_G$. Then there exists a number
  $\lambda=\lambda (M,\phi, d_G)\ge 0$ such that for any initial
  distribution $\mu$ on $E(M)$ we have
\[
\lim_{n\to\infty} \frac{|\phi(w_n)|_G}{n}= \lim_{n\to\infty}
\frac{|\phi(w_n')|_G}{n}=\lambda \text{ almost surely and in } L^1
\text{ with respect to } P_{\mu}.
\]

\end{thm}

\begin{proof}
  Let $\mu_0$ be the unique stationary initial distribution for $M'$.
  As before denote by $\mathcal S:\Omega\to \Omega$ the shift operator
  which erases the first edge of every $\omega=e_1,e_2, \dots \in
  \Omega$.  Stationarity of $\mu_0$ means that $\mathcal S:(\Omega,
  P_{\mu_0})\to (\Omega, P_{\mu_0})$ is a measure-preserving map.
  Since $M'$ is irreducible and aperiodic, $\mathcal S$ is also
  ergodic.

  As before, define $X_n:\Omega\to \mathbb R$ as
\[
X_n(\omega):=|\phi(w_n(\omega))|_G.
\]

Then again it is easy to see that $X_n\ge 0$, $X_{n+m}(\omega)\le
X_n(\omega)+X_m(\mathcal S^n\omega)$. Hence by the Subadditive Ergodic
Theorem there is $\lambda\ge 0$ and there is a subset $Q\subseteq
\Omega$ with $P_{\mu_0}(Z)=0$ such that for any $\omega\not\in Z$

\[
\lim_{n\to\infty} \frac{|\phi(w_n(\omega))|_G}{n}=\lambda.
\]

Let $\mu$ be an arbitrary initial distribution on $E$. Then by
Lemma~\ref{zero} we have $P_{\mu}(Z)=0$.  Thus
\[
\frac{|\phi(w_n)|_G}{n} \to \lambda \text{ almost surely with respect
  to } P_{\mu}.
\]

Note that by the left-invariance of $d_G$ we have $|\phi(w)|_G\le K
|w|$ where $K=\max\{ |\phi(a)|_G : a\in A\}$.  Hence
$X_n/n=|\phi(w_n)|_G/n\le K$ and by the Dominated Convergence Theorem
almost sure convergence of $X_n/n$ implies $L^1$-convergence.

Since $d_G$ is a seminorm on $G$ and the length of any path $w_n'$
differs from $|w_n|$ by at most a fixed constant, it is also true that
$|\phi(w_n')|_G$ differs from $|\phi(w_n)|_G$ by at most a fixed
constant and thus it is also the case that

\[
\lim_{n\to\infty} \frac{|\phi(w_n'(\omega))|_G}{n}=\lim_{n\to\infty}
\frac{|\phi(w_n(\omega))|_G}{n}=\lambda.
\]

\end{proof}

There is substantial flexibility in the choice of the Markov chain
$M'$. The proof of Theorem~\ref{thm:H} goes through without change for
any choice of transition probabilities in $M'$ such that $p_{e,f}>0$
whenever $(e,f)$ is an edge of $\Gamma_1(M)$ and $p_{e,f}=0$ whenever
$(e,f)$ is not an edge of $\Gamma_1(M)$.

\section{Open Problems}\label{sect:prob}

\begin{prob}
  Let $\phi$ be an arbitrary (not necessarily injective) endomorphism
  of $F=F(a_1,\dots, a_k)$. Is $\lambda(\phi)$ rational? Computable?
\end{prob}

\begin{prob}
  Let $\phi\in Aut(F)$. What can be said about the behavior of
  $\lambda(\phi^n)$ as $n\to\infty$? Same for
  $\sqrt[n]{\lambda(\phi^n)}$. How are these quantities connected with
  growth rates of different (or perhaps just top) strata from relative
  train-track representatives of $\phi$?
\end{prob}

It is clear that the asymptotics of $\lambda(\phi^n)$ should reflect
the dynamical properties of $\phi$. For example, it is not hard to see
that for any Nielsen automorphism $\tau$ the stretching factor
$\lambda(\tau^n)$ grows at most linearly and
$\limsup_{n\to\infty}\sqrt[n]{\lambda(\tau^n)}=1$. On the other hand
for hyperbolic automorphisms $\phi$ Theorem~\ref{hyperb} implies that
$\liminf_{n\to\infty}\sqrt[n]{\lambda(\phi^n)}>1$, so that the
sequence $(\lambda(\phi^n))_n$ grows exponentially.

\begin{prob}
  Can one estimate (say in the sense of Large Deviations) the speed of
  convergence $\frac{|\phi(\omega_n)|_G}{n}\to\lambda(\phi)$?

  We have seen that in the case of free group automorphisms for any
  $\epsilon>0$

\[
P_n\big(\frac{|\phi(\omega_n)|}{n}\in (\lambda(\phi)-\epsilon,
\lambda(\phi)+\epsilon)\big)\to 1
\]
with exponentially fast convergence as $n\to\infty$. Are there any
other situations where the speed of convergence in
Theorem~\ref{thm:stretch} can be estimated?
\end{prob}

\begin{prob}
  Let $F=F(a_1,\dots, a_k)$ where $k\ge 2$.  Consider the set
\[
W=\{\lambda(\phi): \phi:F\to Aut(X) \text{ is a free simplicial action
  of } F \text{ on some simplicial tree } X\}.
\]
We know that $W\subseteq \mathbb Q$ and, moreover $2kW\subseteq
\mathbb Z[\frac{1}{2k-1}]$.

Is $W$ a discrete subset of $\mathbb Q$?

\end{prob}

%\begin{prob}
%  Let $F$ be a free group of rank $k>1$ that is realized as the
%  fundamental group of a finite metric graph $X$. Then $\tilde X$ is a
%  metric tree with a free isometric action of $F$ such that $\tilde
%  X/F=X$. Choose $p\in X$ and define a metric $d'$ on $F$ as
%  $d'(f_1,f_2):=d_{\tilde X}(f_1p,f_2p)$.

%  To what extent do the results of this paper carry through to
%  automorphisms of $F$ considered as maps $(F,d')\to (F,d')$? In
%  particular, is it true that in this situation $\lambda(\phi)=1$ if
%  and only if up to a composition with an inner automorphism $\phi$ is
%  induced by an isometry of $X$?
%\end{prob}

\begin{prob}
  The notion of a generic stretching factor for $\phi\in Aut(F)$
  depends on the choice of a free basis $b=(a_1,\dots, a_k)$ of $F$,
  and, more generally, on the choice of a finite generating set $S$ of
  $F$ and the corresponding word metric $d_S$. Denote by
  $\lambda_S(\phi)$ the generic stretching factor of $\phi$ considered
  as a map $(F,d_S)\to (F,d_S)$.

  One can define the following uniform constants
\[
\lambda'(\phi):=\inf\{\lambda_b(\phi): b \text{ is a free basis of }
F\}
\]
and

\[
\lambda''(\phi):=\inf\{\lambda_S(\phi): S \text{ is a finite
  generating set of } F\}.
\]
(Note that $\lambda''(\phi)$ can be defined in the same fashion for an
automorphism $\phi$ of an arbitrary finitely generated group $G$).

For $\phi\in Aut(F)$ are the constants $\lambda'(\phi)$ and
$\lambda''(\phi)$ actually realized by some free bases and finite
generating sets of $F$ accordingly? That is, are the above infima
actually minima? Are $\lambda'(\phi)$ and $\lambda''(\phi)$
algorithmically computable?

Similarly we can define
\[
||\phi||'=\inf\{||\phi||_b: b \text{ is a free basis of } F\}
\]
and
\[
||\phi||''=\inf\{||\phi||_S: S \text{ is a finite generating set of }
F\}.
\]
Since both of these constants are integers, they are clearly
realizable by some $b$ and $S$ accordingly.  Are these constants
algorithmically computable?
\end{prob}

%\begin{prob}
%Let $G$ be a torsion-free hyperbolic group with a word metric $d_G$
%and let $\phi\in Aut(G)$ be such that $\phi: (G, d_G)\to (G,d_G)$ is a
%rough isometry. Does this imply that $[\phi]$ has finite order in
%$Out(G)$? Same question for $\phi$ being a rough similarity.
%\end{prob}

\end{document}